\documentclass[11pt,a4paper]{article}

\usepackage[utf8]{inputenc}
\usepackage[T1]{fontenc}
\usepackage[english]{babel}
\usepackage{amsmath}
\usepackage{amsfonts}
\usepackage{amsthm}
\usepackage{amssymb}
\usepackage{makeidx}
\usepackage{graphicx}
\usepackage{lmodern}
\usepackage[left=2cm,right=2cm,top=2cm,bottom=2cm]{geometry}

\usepackage{color}
\usepackage{comment}
\usepackage[dvipsnames]{xcolor}
\usepackage{graphicx}
\usepackage{framed}
\usepackage{tikz}
\usepackage{calrsfs}
\DeclareMathAlphabet{\pazocal}{OMS}{zplm}{m}{n}

\usepackage{fourier-orns}
\usepackage{mathtools}
\usepackage{relsize}
\usepackage{dsfont}
\usepackage{comment}
\usepackage{hyperref}

\newcommand{\R}{\mathbb{R}}

\renewcommand{\H}{\mathbb{H}}
\newcommand{\J}{\mathbb{J}}

\newcommand{\F}{\pazocal{F}}
\newcommand{\G}{\pazocal{G}}
\newcommand{\U}{\pazocal{U}}
\newcommand{\K}{\pazocal{K}}
\newcommand{\V}{\pazocal{V}}

\newcommand{\Wcal}{\mathcal{W}}
\newcommand{\Lcal}{\mathcal{L}}
\newcommand{\Vcal}{\mathcal{V}}
\newcommand{\Pcal}{\mathcal{P}}

\newcommand{\Id}{\textnormal{Id}}
\newcommand{\D}{\textnormal{D}}

\newcommand{\Var}{\textnormal{Var}}
\newcommand{\supp}{\textnormal{supp}}

\newcommand{\Lip}{\textnormal{Lip}}

\newcommand{\Bgamma}{\boldsymbol{\gamma}}

\newcommand{\Bmu}{\boldsymbol{\mu}}
\newcommand{\Bpartial}{\boldsymbol{\partial}}

\newcommand{\INTDom}[3]{\int_{#2} #1 \textnormal{d} #3}
\newcommand{\INTSeg}[4]{\int_{#3}^{#4} #1 \textnormal{d} #2}
\newcommand{\NormL}[3]{\parallel \hspace{-0.1cm} #1 \hspace{-0.1cm} \parallel _ {L^{#2}(#3)}}
\newcommand{\NormC}[3]{\left\| #1  \right\| _ {C^{#2}(#3)}}

\newcommand{\derv}[2]{\frac{\textnormal{d} #1}{ \textnormal{d} #2}}

\newcommand{\BGamma}{\mbox{$ \textnormal{l} \hspace{-0.1em}\Gamma\hspace{-0.88em}\raisebox{-0.98ex}{\scalebox{2}
  {$\color{white}.$}}\hspace{0.46em}$}{}}  
  
\newtheorem{rmk}{Remark}
\newtheorem{lem}{Lemma}
\newtheorem{Def}{Definition}
\newtheorem{thm}{Theorem}
\newtheorem{prop}{Proposition}

\newtheorem{cor}{Corollary}

\title{The Pontryagin Maximum Principle in the Wasserstein Space\footnote{This work has been carried out in the framework of Archimède
Labex (ANR-11-LABX-0033) and of the A*MIDEX project (ANR-
11-IDEX-0001-02), funded by the "Investissements d'Avenir" French
Government program managed by the French National Research
Agency (ANR). The authors also thank the reviewer for the several useful comments that he/she provided.}}
\author{Benoît Bonnet\footnote{Aix Marseille Universit\'e, CNRS, ENSAM, Universit\'e de Toulon, LSIS, Marseille, France} , Francesco Rossi\footnote{Dipartimento di Matematica "Tullio Levi-Civita" Universit\`a degli Studi di Padova, Padova,
Italy}}

\date{January 18, 2018} 

\begin{document}
\maketitle

\begin{abstract}
We prove a Pontryagin Maximum Principle for optimal control problems in the space of probability measures, where the dynamics is given by a transport equation with non-local velocity. We formulate this first-order optimality condition using the formalism of subdifferential calculus in Wasserstein spaces. We show that the geometric approach based on needle variations and on the evolution of the covector (here replaced by the evolution of a mesure on the dual space) can be translated into this formalism.
\end{abstract}

\section{Introduction}
\label{section:Introduction}

Transport equations with non-local interaction terms have been intensively studied for decades by various communities. They were for instance already introduced in statistical physics in 1938 when Vlasov proposed these equations to describe long-range Coulomb interactions \cite{vlasov}. For such reasons, several transport equations appear as mean-field limits of particle systems, see e.g. \cite{golse,spohn}. More recently, the study of crowd modelling has stimulated a renewed interest for these equations. Indeed, pedestrians have a long-range perception of their space, and thus choose their path based on long-range interactions. While such interactions do not enjoy action-reaction properties which are typical in physical models, methods connected to mean-field limit approaches have shown their adaptability in this setting too (see e.g. \cite{Bellomo2013,CPT,mauryvenel1,mauryvenel2}). More generally, the study of other kind of interacting agents, such as opinion dynamics on networks \cite{Bellomo2013,HK}, or animal flocks \cite{ballerini,CS2}, has been conveyed with similar techniques.

Several contributions have shown that the natural setting for studying transport equations with non-local terms is the space of measures endowed with the Wasserstein distance, see e.g. \cite{AGS}. In this case, existence and uniqueness of the solution of a Cauchy problem are ensured by a natural Lipschitz condition \cite{AmbrosioGangbo}, and metric estimates for the associated flow are available too \cite{Pedestrian}. For simplicity, we will only deal with measures with compact support, for which the Wasserstein distance is always finite.

Beside the analysis of such partial differential equations, it is now of great interest to study {\bf control problems for the transport equation with non-local velocities}. Apart from a few recent results about controllability \cite{Duprez2019}, most of the contributions in this direction have considered optimal control problems, i.e. the minimization of a functional where the constraint is a controlled dynamics. Applications of these problems are of great interest, canonical examples being provided e.g. by the minimal escape time problem for a crowd \cite{albicristiani,lifebelt} or the enforcement of consensus in a network by minimizing the variance of the opinions (see e.g. \cite{SparseJQMF,SparseJQ,ControlKCS}).

Existence of optimal controls has been investigated in \cite{MFOC}, as well as in the setting of mean-field control \cite{achdou2,achdou1}. Convergence of optimizers via the mean-field limit of the dynamics was also studied with methods related to $\Gamma$-convergence in \cite{FPR}.

The next logical step in this study is the derivation of first-order necessary optimality conditions allowing to characterize and compute optimal trajectories. Although Hamilton-Jacobi optimality conditions in Wasserstein spaces have received some attention, see e.g. the seminal paper \cite{HJBWasserstein} and recent developments in the field of control theory \cite{Cavagnari2018}, Pontryagin optimality conditions remain rather unexplored. A first result in this direction was presented in \cite{MFPMP}, in which a coupled PDE-ODE system was studied, in which the control acts on the ODE part only. The main result was a necessary first-order condition written as a Pontryagin Maximum Principle. Instead, we turn our attention here to a control problem formulated directly on the PDE. Then, one needs a sufficiently rich differential structure to compute derivatives of the functional to be minimized with respect to the control. In this context, the state is represented by a measure, for which the adapted setting is given by subdifferential calculus in Wasserstein spaces. We recall the main useful results of this theory in Section \ref{s-wass} (see also \cite{AGS} for a thorough introduction).

Our contribution in this article is to show that, in this general framework, several results of geometric control can be translated from finite-dimensional dynamical systems to transport equations with non-local velocities. With this aim, we derive a new Pontryagin Maximum Principle in this infinite-dimensional setting. While the proof scheme is close to the classical finite-dimensional case, each step requires the definition of tools adapted to Wasserstein spaces and additional technical care in the different arguments.

As a result, the new Pontryagin Maximum Principle (PMP in the following) is formulated in the language of subdifferential calculus in Wasserstein spaces.  In particular, the state-costate variables are here replaced by a measure on the product of the tangent and cotangent bundle. The dynamics is given by an Hamiltonian system in the space of measures, similar to what studied in \cite{AmbrosioGangbo}, where the corresponding Hamiltonian is given by a maximization in an adapted space of controls functions satisfying Lipschitz constraints.\\

In the sequel, we shall study Pontryagin-type optimality condition for optimal control problems given in the general form
\begin{equation} \label{eq:General_OCP}
(\pazocal{P}) ~~ \left\{
\begin{aligned}
& \underset{u \in \U}{\text{min}} \left[ \INTSeg{L(\mu(t),u(t))}{t}{0}{T} +  \varphi(\mu(T)) \right], \\
& \text{s.t.} 
\left\{ 
\begin{aligned}
& \partial_t \mu(t) + \nabla \cdot ((v[\mu(t)](t,\cdot) + u(t,\cdot)) \mu(t)) = 0, \\
& \mu(0) = \mu^0 \in \Pcal_c(\R^d).
\end{aligned}
\right.
\end{aligned}
\right.
\end{equation}
As already stated, our formulation of the PMP deeply relies on the formalism of subdifferential calculus in Wasserstein spaces (see e.g. \cite{CarrilloLisiniMainini,GianazzaSavareToscani}). In this formalism, the \textit{extended subdifferential} $\Bpartial \phi(\mu)$ (see Definition \ref{def:Subdifferentials} below) of a functional $\phi(\cdot)$ at a given measure $\mu \in \Pcal_c(\R^d)$ is made of transport plans. As it is the case for subdifferential calculus in Banach spaces, there exists a notion of \textit{minimal selection} (see Theorem \ref{thm:subdifferential_localslope} below) among the elements of this subdifferential. The minimal selection in an extended subdifferential, which we denote by $\Bpartial^{\circ} \phi(\mu)$, plays the same conceptual role of the \textit{gradient} of a differentiable functional. The existence of such minimal selection is a consequence of the \textit{regularity} hypothesis (see Definition \ref{def:Regular} and the corresponding Theorem \ref{thm:subdifferential_localslope} below), that we impose to the functionals studied in the following. 

In this context, the \textit{barycenter} $\bar{\Bgamma}_{\phi}^{\circ} : \R^d \rightarrow \R^d$ (see Definition \ref{def:Barycenter} below) of the minimal selection is the closest object to what would be a gradient in the sense of subdifferential calculus, in particular when computing derivatives along curves of measures (see Proposition \ref{prop:Chainrule} below). However, barycenters of extended subdifferentials are not in the  \textit{classical subdifferentials} in general. Yet for a good score of functionals involved in applications such as potential and interaction energies, relative entropies, variance functionals (see e.g. Section \ref{section:Examples} below for some examples), the minimal selection is induced by its barycenter. In this case, the latter is referred to as the \textit{Wasserstein gradient} (see Definition \ref{def:Wass_Grad} below) $\nabla_{\mu}\phi(\mu) : \R^d \rightarrow \R^d$ of the functional $\phi(\cdot)$ at $\mu$.\\

We introduce in Theorem \ref{thm:Introduction_PMP} below a heuristic version of our main result and postpone for the sake of readability its precise statement to Section \ref{section:PMP}, Theorem \ref{thm:General_PMP}. In the sequel, we will denote by $B_{2d}(0,R)$ the ball of radius $R$ centered at $0$ in $\R^{2d}$, by $\pi^1,\pi^2 : \R^{2d} \rightarrow \R^d$ the projection operators on the first and second components and by $K$ a generic compact subset of $\R^d$. 

\begin{thm}[Heuristic statement of the Pontryagin Maximum Principle for $(\pazocal{P})$] 
\label{thm:Introduction_PMP}
Let $(u^*(\cdot),\mu^*(\cdot)) \in \U \times \Lip([0,T],\Pcal_c(\R^d))$ be an optimal pair control-trajectory for $(\pazocal{P})$ and assume that hypotheses \textnormal{\textbf{(H)}} of Theorem \ref{thm:General_PMP} below hold.

Then, there exist a constant $R > 0$ depending on $\mu^0$, $T$, $\U$, $v[\cdot](\cdot,\cdot)$, $\varphi(\cdot)$, $L(\cdot,\cdot)$ and  a curve $\nu^*(\cdot) \in \Lip([0,T],\Pcal(\overline{B_{2d}(0,R)})$ Lipschitzian with respect to the $W_1$-metric satisfying the following conditions : 
\begin{enumerate}
\item[(i)] It solves the forward-backward system of continuity equations
\begin{equation*}
\left\{
\begin{aligned}
& \partial_t \nu^*(t) + \nabla_{(x,r)} \cdot \left( \J_{2d} \tilde{\nabla}_{\nu} \H_c(t,\nu^*(t),u^*(t))\nu^*(t) \right) = 0 ~ \text{in $[0,T] \times \R^{2d}$}, \\
& \pi^1_{\#} \nu^*(0) = \mu^0, \\
& \pi^2_{\#} \nu^*(T) = (-\bar{\Bgamma}^{\circ}_{\varphi})_{\#} \mu^*(T),
\end{aligned}
\right.
\end{equation*}
where the vector field $\tilde{\nabla}_{\nu} \H_c(t,\nu^*(t),u^*(t))(\cdot,\cdot)$ is (almost) the Wasserstein gradient of a suitable compactification of the  \textnormal{infinite dimensional Hamiltonian} $\H(\cdot,\cdot,\cdot)$ of the system, defined by 
\begin{equation*}
\H(t,\nu,\omega) = \INTDom{\langle r , v[\pi^1_{\#}\nu](t,x) + \omega(x) \rangle}{\R^{2d}}{\nu(x,r)} - L(\pi^1_{\#}\nu,\omega)
\end{equation*}
for any $(t,\nu,\omega) \in [0,T] \times \Pcal_c(\R^{2d}) \times U$.

\item[(ii)] It satisfies the Pontryagin maximization condition
\begin{equation*}
\hspace{-0.6cm} \H_c(t,\nu^*(t),u^*(t)) = \max\limits_{\omega \in U} \left[ \H_c(t,\nu^*(t),\omega) \right]
\end{equation*}
for $\Lcal^1$-almost every $t \in [0,T]$.
\end{enumerate}
\end{thm}

The structure of the article is the following : in Section \ref{section:WassersteinRecall} we recall useful results of analysis in Wasserstein spaces, PDEs with non-local velocities and subdifferential calculus in $(\Pcal_2(\R^d),W_2)$. We also prove in Proposition \ref{prop:Directional_derivative_flow} an existence and characterization result for directional derivatives along measure curves for non-local flows. In Section \ref{section:PMP} we state and prove our main result. We first introduce in Section \ref{subsection:SimplePMP} the main steps of our proof strategy - in particular the concept of \textbf{needle like variation} -, on a simpler instance $(\pazocal{P}_1)$ of problem $(\pazocal{P})$. We proceed to prove Theorem \ref{thm:General_PMP} in Section \ref{subsection:GeneralPMP}. In Section \ref{section:Examples} we discuss more in details the set of hypotheses \textbf{(H)} of Theorem \ref{thm:General_PMP} and list some relevant examples of classical functionals satisfying them.

%%%%%%%%%%%%%%%%%%%%%%%%%%%%%%%%%%%%%%%%%%%%%%%%%%%%%%%%%%%%%%%%%%%%%%%%%%%%
%
%									NEW SECTION
%
%%%%%%%%%%%%%%%%%%%%%%%%%%%%%%%%%%%%%%%%%%%%%%%%%%%%%%%%%%%%%%%%%%%%%%%%%%%%

\section{Analysis in Wasserstein spaces} \label{s-wass}
\label{section:WassersteinRecall}

In this section, we recall several notions about analysis in the space of probability measures, optimal transport theory, Wasserstein spaces, continuity equations and subdifferential calculus in the space $(\Pcal_2(\R^d),W_2)$. All the results stated in this section are well-known, at the exception of Proposition  \ref{prop:Directional_derivative_flow} which is a generalization of the classical differentiation result for smooth flows of diffeomorphisms that we recall in Proposition \ref{prop:Differential_local_flow}.

%%%%%%%%%%%%%%%%%%%%%%%%%%%%%%%%%%%%%%%%%%%%%%%%%%%%%%%%%%%%%%%%%%%%%%%%%%%%

\subsection{The optimal transport problem and Wasserstein spaces}

In this section, we introduce some classical notations and results of optimal transport and analysis in Wasserstein spaces.

We denote by $\Pcal(\R^d)$ the space of Borel probability measures over $\R^d$ and by $\Lcal^d$ the standard Lebesgue measure on $\R^d$. For $p \geq 1$, we define $\Pcal_p(\R^d)$ as the subset of $\Pcal(\R^d)$ of measures having finite $p$-th moment, i.e.
\begin{equation*}
\Pcal_p(\R^d) = \left\{ \mu \in \Pcal(\R^d) ~\text{s.t.}~ \INTDom{|x|^p}{\R^d}{\mu(x)} < +\infty \right\}.
\end{equation*}

The \textit{support} of a Borel probability measure $\mu \in \Pcal(\R^d)$ is defined as the closed set $\supp(\mu) = \{ x \in \R^d ~\text{s.t.}~ \mu(\pazocal{N}) > 0 ~\text{for any neighbourhood $\pazocal{N}$ of $x$}\}$. We denote by $\Pcal_c(\R^d)$ the subset of $\Pcal(\R^d)$ of measures which supports are compact. 

We say that a sequence $(\mu_n) \subset \Pcal(\R^d)$ of Borel probability measures \textit{converges narrowly} towards $\mu \in \Pcal(\R^d)$, denoted by $\mu_n \underset{n \rightarrow +\infty}{\rightharpoonup} \mu$, provided that
\begin{equation}
\label{eq:Narrow_convergence}
\INTDom{\phi(x)}{\R^d}{\mu_n(x)} ~\underset{n \rightarrow + \infty}{\longrightarrow}~ \INTDom{\phi(x)}{\R^d}{\mu(x)} ~~ \text{for all $\phi \in C^0_b(\R^d)$}
\end{equation}
where $C^0_b(\R^d)$ denotes the set of continuous and bounded functions from $\R^d$ into $\R$. 

We recall the definitions of \textit{pushforward} of a Borel probability measure through a Borel map and \textit{transport plan}. 

\begin{Def}[Pushforward of a measure through a Borel map]
\hfill \\ Given a Borel probability measure $\mu \in \Pcal(\R^d)$ and a Borel map $f : \R^d \rightarrow \R^d$, the \textit{pushforward} $f_{\#} \mu$ of $\mu$ through $f(\cdot)$ is defined as the only Borel probability measure such that $f_{\#} \mu (B) = \mu(f^{-1}(B))$ for any Borel set $B \subset \R^d$. 
\end{Def}

\begin{Def}[Transport plan]
Given two probability measures $\mu$ and $\nu$ on $\R^d$, we say that $\gamma \in \Pcal(\R^{2d})$ is a \textit{transport plan} between $\mu$ and $\nu$, denoted by $\gamma \in \Gamma(\mu,\nu)$, provided that $\gamma(A \times \R^d) = \mu(A)$ and $\gamma(\R^d \times B) = \nu(B)$ for any Borel subsets $A,B \subset \R^d$, or equivalently $\pi^1_{\#} \gamma = \mu$ and $\pi^2_{\#} \gamma = \nu$.

Given a probability measure $\gamma \in \R^{2d}$, we also denote by $\Gamma(\gamma,\nu)$ the set of plans $\Bmu \in \Pcal(\R^{3d})$ such that $\pi^{1,2}_{\#} \Bmu = \gamma$ and $\pi^3_{\#} \Bmu = \nu$ where $\pi^{1,2} : (x,y,z) \in \R^{3d} \mapsto (x,y) \in \R^{2d}$.
\end{Def}

We recall in the following Proposition three useful convergence results for sequences of probability measures and functions (see e.g. \cite[Chapter 5]{AGS}).

\begin{prop}[Convergence results]
\label{prop:Convergence_results}
Let $(\mu_n) \subset \Pcal(\R^d)$ be a sequence narrowly converging to $\mu \in \Pcal(\R^d)$, $(f_n)$ be a sequence of $\mu$-measurable functions pointwisely converging to $f$ and $g \in C^0(\R^d)$. 
\begin{enumerate}
\item[(i)] Suppose that $x \mapsto |g(x)|$ is \textnormal{uniformly integrable} with respect to the family $\{ \mu_n \}_{n=1}^{\infty}$, i.e 
\begin{equation*}
\lim\limits_{k \rightarrow +\infty} \INTDom{|g(x)|}{ \{ |g(x)| \geq k \} }{\mu_n(x)} = 0
\end{equation*}
for all $n \geq 1$. Then, the sequence $(\INTDom{g(x)}{\R^d}{\mu_n(x)}) \subset \R$ converges to $\INTDom{g(x)}{\R^d}{\mu(x)}$ as $n \rightarrow +\infty$. \\
\item[(ii)] The sequence $(g_{\#}\mu_n) \subset \Pcal(\R^d)$ narrowly converges to $g_{\#}\mu$ as $n \rightarrow +\infty$. \\
\item[(iii)] \textnormal{(Vitali convergence theorem)} Suppose that the family $x \mapsto |f_n(x)|$ is \textnormal{uniformly integrable} with respect to the measure $\mu$, i.e. 
\begin{equation*}
\lim\limits_{k \rightarrow +\infty} \INTDom{|f_n(x)|}{ \{ |f_n(x)| \geq k \} }{\mu(x)} = 0
\end{equation*}
for all $n \geq 1$ and also assume that $|f(x)| < +\infty$ for $\mu$-almost every $x \in \R^d$. Then $(f_n)$ converges uniformly to $f$ in $L^1(\R^d;\mu)$ as $n \rightarrow +\infty$.
\end{enumerate}
\end{prop}

In the 40's, Kantorovich introduced the \textit{optimal mass transportation problem} in its modern mathematical formulation : given two probability measures $\mu,\nu \in \Pcal(\R^d)$ and a cost function $c : \R^{2d} \rightarrow \R$, find a \textit{transport plan} $\gamma \in \Gamma(\mu,\nu)$ such that 
\begin{equation*}
\INTDom{c(x,y)}{\R^{2d}}{\gamma(x,y)} = \min \left\{ \INTDom{c(x,y)}{\R^{2d}}{\gamma'(x,y)} ~~ \text{s.t.} ~ \gamma' \in \Gamma(\mu,\nu) \right\}.
\end{equation*}

This problem has been extensively studied in very broad contexts (see e.g. \cite{AGS,villani1}) with high levels of generality on the underlying spaces and cost functions.  In the particular case where $c(x,y) = |x-y|^p$ for some real number $p \geq 1$, the optimal transport problem can be used to define a distance over the subspace $\Pcal_p(\R^d)$ of $\Pcal(\R^d)$.

\begin{Def}[Wasserstein distance and Wasserstein spaces]
\hfill \\ Given two probability measures $\mu,\nu \in \Pcal_p(\R^d)$, the $p$-Wasserstein distance $W_p$ between $\mu$ and $\nu$ is defined by
\begin{equation*}
W_p(\mu,\nu) = \min \left\{ \left( \INTDom{|x-y|^p}{\R^{2d}}{\gamma(x,y)} \right)^{1/p} ~~ \text{s.t.} ~ \gamma \in \Gamma(\mu,\nu) \right\}.
\end{equation*}

The set of plans $\gamma \in \Gamma(\mu,\nu)$ achieving this optimal value is denoted\footnote{We omit the dependence on $p$ for clarity and conciseness.} by $\Gamma_o(\mu,\nu)$ and referred to as the set of \textnormal{optimal transport plans} between $\mu$ and $\nu$. The space $(\Pcal_p(\R^d),W_p)$ of probability measures with finite $p$-th moment endowed with the $p$-th Wasserstein metric is called the \textnormal{Wasserstein space} of order $p$.
\end{Def}

We recall some of the interesting properties of these spaces in the following Proposition (see e.g. \cite[Chapter 7]{AGS} or \cite[Chapter 6]{villani1}).

\begin{prop}[Properties of the Wasserstein distance]
\label{prop:Properties_Wp}
The topology induced in $\Pcal_p(\R^d)$ by the Wasserstein metric $W_p$ metrizes the weak-* topology of probability measures induced by the narrow convergence \eqref{eq:Narrow_convergence}. More precisely,
\begin{equation*}
W_p(\mu_n,\mu) \underset{n \rightarrow +\infty}{\longrightarrow} 0 ~\Longleftrightarrow ~ \mu_n \underset{n \rightarrow +\infty}{\rightharpoonup} \mu ~\text{and}~ \INTDom{|x|^p}{\R^d}{\mu_n(x)} \underset{n \rightarrow +\infty}{\longrightarrow} \INTDom{|x|^p}{\R^d}{\mu(x)}
\end{equation*}

For compactly supported measures $\mu,\nu \in \Pcal_c(\R^d)$, the Wasserstein distances are ordered, i.e. $p_1 \leq p_2 ~\Longrightarrow~ W_{p_1}(\mu,\nu) \leq W_{p_2}(\mu,\nu)$. In particular when $p = 1$, the following \textnormal{Kantorovich-Rubinstein duality formula} holds 
\begin{equation} \label{eq:Kantorovich_duality}
W_1(\mu,\nu) = \sup \left\{ \INTDom{\phi(x) \,}{\R^d}{(\mu-\nu)(x)} ~\text{s.t.}~ \Lip(\phi,\R^d) \leq 1 ~\right\}.
\end{equation}
\end{prop}

In what follows, we shall mainly restrict our considerations to the Wasserstein spaces of order 1 and 2 built over $\Pcal_c(\R^d)$. We end this introductory paragraphs by recalling the concepts of \textit{disintegration} and \textit{barycenter} in the context of optimal transport. 

\begin{Def}[Disintegration and barycenter]
\label{def:Barycenter}
Let $\mu,\nu \in \Pcal_p(\R^d)$ and $\gamma \in \Gamma(\mu,\nu)$ be a transport plan between $\mu$ and $\nu$. We define the \textit{disintegration} $\{ \gamma_x\}_{x \in \R^d} \subset \Pcal_p(\R^d)$ of $\gamma$ on its first marginal $\mu$, usually denoted by $\gamma = \int \gamma_x \textnormal{d} \mu(x)$, as the $\mu$-almost uniquely determined Borel family of probability measures such that 
\begin{equation*}
\INTDom{\phi(x,y)}{\R^{2d}}{\gamma(x,y)} = \INTDom{\INTDom{\phi(x,y)}{\R^d}{\gamma_x(y)}}{\R^d}{\mu(x)}, 
\end{equation*}
for any Borel map $\phi : \R^{2d} \rightarrow \R^d$. 

The \textit{barycenter} $\bar{\gamma} \in L^p(\R^d,\R^d;\mu)$ of the plan $\gamma$ is then defined by 
\begin{equation*}
\bar{\gamma} : x \in \supp(\mu) \mapsto \INTDom{y \,}{\R^d}{\gamma_x(y)}.
\end{equation*}
\end{Def}

%%%%%%%%%%%%%%%%%%%%%%%%%%%%%%%%%%%%%%%%%%%%%%%%%%%%%%%%%%%%%%%%%%%%%%%%%%%%

\subsection{The continuity equation with non-local velocities on $\R^d$}

In this section, we introduce the continuity equations with non-local velocities in $(\Pcal_c(\R^d),W_1)$. These equations write
\begin{equation} \label{eq:Nonlocal_continuity_equation}
\partial_t \mu(t) + \nabla \cdot \left( v[\mu(t)](t,\cdot) \mu(t) \right) = 0,
\end{equation}
where $t \mapsto \mu(t)$ is a narrowly continuous family of probability measures on $\R^d$ and $(t,x) \mapsto v[\mu](t,x)$ is a Borel family of vector fields for any $\mu \in \Pcal_c(\R^d)$, satisfying the condition
\begin{equation}
\INTSeg{\INTDom{|v[\mu(t)](t,x)| \,}{\R^d}{\mu(t)(x)}}{t}{0}{T} < +\infty.
\end{equation}  
Equation \eqref{eq:Nonlocal_continuity_equation} has to be understood in the sense of distributions, i.e.  
\begin{equation}
\label{eq:NonlocalPDE_distributions1}
\INTSeg{\INTDom{\left( \partial_t \phi(t,x) + \left\langle \nabla_x \phi(t,x) , v[\mu(t)](t,x) \right\rangle \right)}{\R^d}{\mu(t)(x)}}{t}{0}{T} = 0
\end{equation}
for all $\phi \in C^{\infty}_c([0,T] \times \R^d)$, or alternatively as
\begin{equation}
\label{eq:NonlocalPDE_distributions2}
\derv{}{t} \INTDom{\phi(x)}{\R^d}{\mu(t)(x)} = \INTDom{\langle \nabla \phi(x) , v[\mu(t)](t,x) \rangle}{\R^d}{\mu(t)(x)}
\end{equation}
for all $\phi \in C^{\infty}_c(\R^d)$ and $\Lcal^1$-almost every $t \in [0,T]$. 

As already mentioned in the introduction, these equations are interesting for a large number of applications. It is important to notice that $v[\mu]$ depends on the whole measure $\mu$ and not only on its values at some points as it is usually the case for non-linear conservation laws.

We now recall a theorem which was first derived in \cite{AmbrosioGangbo} providing existence, uniqueness and representation formula for solutions of \eqref{eq:Nonlocal_continuity_equation}. We state here a version explored in  \cite{Pedestrian,ControlKCS} that is more suited to our control-theoretic framework.

\begin{thm}[Existence, uniqueness and representation of solutions for \eqref{eq:Nonlocal_continuity_equation}] \label{thm:Existence_uniqueness_non-local_PDE}
Consider a non-local velocity field $v[\cdot](\cdot,\cdot)$ defined as
\begin{equation}
v : \mu \in \Pcal_c(\R^d) \mapsto v[\mu](\cdot,\cdot) \in L^{\infty}(\R , C^1 \cap L^{\infty} (\R^d,\R^d)), 
\end{equation}
and satisfying the following assumptions
\begin{framed}
\vspace{-0.3cm}
\begin{center}
\textnormal{\textbf{(H')}}
\end{center}
\begin{enumerate}
\item[$\diamond$] There exists positive constants $L_1$ and $M$ such that
\begin{equation*}
| v[\mu](t,x) - v[\mu](t,y) | \leq L_1 |x-y| ~~ \text{and} ~~ | v[\mu](t,x) | \leq M (1+ |x|)
\end{equation*}
for every $\mu \in \Pcal_c(\R^d)$, $t \in \R$ and $(x,y) \in \R^{2d}$;
\item[$\diamond$] There exists a positive constant $L_2$ such that
\begin{equation*}
\NormC{v[\mu](t,\cdot) - v[\nu](t,\cdot)}{0}{\R^d} \leq L_2 W_1(\mu,\nu) 
\end{equation*}
for every $\mu,\nu \in \Pcal_c(\R^d)$ and $t \in \R$;
\end{enumerate}
\vspace{-0.3cm}
\end{framed}

Then for every initial datum $\mu^0 \in \Pcal_c(\R^d)$, the Cauchy problem
\begin{equation} \label{eq:Nonlocal_Cauchy_problem}
\left\{
\begin{aligned}
& \partial_t \mu(t) + \nabla \cdot \left( v[\mu(t)](t,\cdot) \mu(t) \right) = 0 \\
& \mu(0) = \mu^0,
\end{aligned}
\right.
\end{equation}
admits a unique solution $\mu(\cdot)$ in $C^0(\R,\Pcal_c(\R^d))$. This solution is locally Lipschitz in $t$ with respect to the $W_1$-metric. Besides, if $\mu^0$ is absolutely continuous with respect to $\Lcal^d$, then $\mu(t)$ is absolutely continuous with respect to $\Lcal^d$ as well for all times $t \geq 0$. \\
Furthermore for every $T > 0$ and every $\mu^0,\nu^0 \in \Pcal_c(\R^d)$, there exists $R_T > 0$ depending on $\supp(\mu^0)$ and $C_T > 0$ such that
\begin{equation*}
\supp(\mu(t)) \subset \overline{B(0,R_T)} ~~ \text{and} ~~ W_1(\mu(t),\nu(t)) \leq C_T W_1(\mu^0,\nu^0),
\end{equation*}
for all times $t \in [0,T]$ and any solutions $\mu(\cdot),\nu(\cdot)$ of \eqref{eq:Nonlocal_Cauchy_problem}.

Let $(\Phi_{(0,t)}^v[\mu^0](\cdot))_{t \geq 0}$ be the family of flows of diffeomorphisms generated by the non-local vector field $v[\mu(t)](t,\cdot)$, defined as the unique solution of 
\begin{equation}
\label{eq:Flow_def}
\left\{
\begin{aligned}
\partial_t \Phi_{(0,t)}^v[\mu^0](x) & = v[\mu(t)] \left( t,\Phi_{(0,t)}^v[\mu^0](x) \right), \\
\Phi^v_{(0,0)}[\mu^0](x) & = x \hspace{1.25cm} \text{for all $x$ in $\R^d$}.
\end{aligned}
\right.
\end{equation}
Then, the unique solution of the Cauchy problem \eqref{eq:Nonlocal_Cauchy_problem} can be expressed at time $t$ as $\mu(t) = \Phi_{(0,t)}^v[\mu^0](\cdot)_{\#} \mu^0$. 
\end{thm}

We recall below a standard result which links the differential of the flow of diffeomorphisms of an ODE at time $t$ to the solution of a corresponding linearized Cauchy problem (see e.g. \cite{BressanPiccoli}). 

\begin{prop}[Differential of a flow]
\label{prop:Differential_local_flow}
Let $(t,x) \mapsto v(t,x)$ be measurable in $t$ as well as sublinear and $C^1$ in $x$. Define the family of $C^1$-flows $(\Phi^v_t(\cdot))_{t \geq 0}$ associated to $v(\cdot,\cdot)$ by \eqref{eq:Flow_def} in the case where $v(\cdot,\cdot)$ is independent from $\mu(\cdot)$. 

Then, it holds that the differential $\D_x \Phi^v_{(s,t)}(x)\cdot h$ of the flow between times $s$ and $t$, evaluated at $x$ and applied to some vector $h \in \R^d$ is the unique solution $w(\cdot,x)$ of the linearized Cauchy problem
\begin{equation*}
\partial_t w(t,x) = \D_x v(t,\Phi^v_{(s,t)}(x)) \cdot w(t,x) ~,~ w(s,x) = h. 
\end{equation*}
\end{prop}

This characterization is essential for proving the Pontryagin Maximum Principle in the usual finite dimensional setting using the needle-like variations approach. We shall prove in Proposition \ref{prop:Directional_derivative_flow} a generalization of this result in the non-local case where the initial measure is perturbed by a Lipschitz family of continuous and bounded maps. Such a result is crucial to study the first order perturbation induced by a needle-like variation on a measure curve in the non-local setting.

\subsection{Subdifferential calculus in $(\Pcal_2(\R^d),W_2)$}

In this section, we recall some elements of subdifferential calculus in the Wasserstein space $(\Pcal_2(\R^d),W_2)$. For a thorough introduction, see \cite[Chapters 9-11]{AGS} where the full theory is developed and applied to the study of gradient flows. 

Throughout this section, we denote by $\phi : \Pcal_2(\R^d) \rightarrow (-\infty,+\infty]$ a proper, lower-semicontinuous functional. We denote the effective domain $D(\phi)$ of $\phi(\cdot)$ as the set of points where it is finite, i.e.
\begin{equation*}
D(\phi) = \{ \mu \in \Pcal_2(\R^d) ~\text{s.t.}~ \phi(\mu) < + \infty\}. 
\end{equation*}

We further assume that for $\tau_* > 0$ small enough, the \textit{Moreau-Yosida} relaxation of $\phi(\cdot)$ defined by
\begin{equation} \label{eq:Moreau_Yosida}
\phi_{\text{M}}(\mu,\tau;\cdot) : \nu \mapsto \frac{1}{2\tau} W_2^2(\mu,\nu) + \phi(\nu)
\end{equation}
attains a minimum at some $\mu_{\tau} \in D(\phi)$ for any $\tau \in (0,\tau_*)$. This technical assumption is satisfied whenever $\phi(\cdot)$ is bounded from below and at least lower-semicontinuous and is crucial for proving the main results of the theory developed in \cite[Chapter 10]{AGS}. 

We start by introducing the concept of \textit{extended subdifferentials} for a functional defined over the Wasserstein space $(\Pcal_2(\R^d),W_2)$. 

\begin{Def}[Extended subdifferential]
\label{def:Subdifferentials}
Let $\mu^1 \in D(\phi)$. We say that a transport plan $\Bgamma \in \Pcal_2(\R^{2d})$ belongs to the \textit{extended (Fr\'echet) subdifferential} $\Bpartial \phi(\mu^1)$ of $\phi(\cdot)$ at $\mu^1$ provided that 
\begin{enumerate}
\item[(i)] $\pi^1_{\#} \Bgamma = \mu^1$, 
\item[(ii)] for all  $\mu^3 \in \Pcal_2(\R^d)$ it holds that 
\begin{equation*}
\phi(\mu^3) - \phi(\mu^1) ~ \geq ~ \inf\limits_{\Bmu \in \Gamma^{1,3}_o(\Bgamma,\mu^3)} \left[ \INTDom{\langle x_2 , x_3-x_1 \rangle}{\R^{3d}}{\Bmu} \right] + o(W_2(\mu^1,\mu^3)),
\end{equation*}
where $\Gamma^{1,3}_o(\Bgamma,\mu^3) = \{\Bmu \in \Gamma(\Bgamma,\mu^3) ~\text{s.t.}~ \pi^{1,3}_{\#} \Bmu \in \Gamma_o(\mu^1,\mu^3) \}$. 
\end{enumerate}
Moreover, we say that an extended subdifferential $\Bgamma$ is \textit{induced by a plan} if there exists $\xi \in L^2(\R^d,\R^d;\mu^1)$ such that $\Bgamma = (\Id \times \xi)_{\#} \mu^1$. In which case, $\xi(\cdot)$ belongs to the \textit{classical subdifferential} $\partial \phi(\mu)$ of $\phi(\cdot)$ at $\mu$.   

We say that a transport plan $\Bgamma \in \Pcal_2(\R^{2d})$ belongs to the \textit{strong extended subdifferential} $\Bpartial_S \phi(\mu^1)$ of $\phi(\cdot)$ at $\mu^1$ if the following stronger condition holds:
\begin{equation} 
\label{eq:Def_StrongSubdiff} 
\phi(\mu^3) - \phi(\mu^1) ~ \geq ~ \INTDom{\langle x_2 , x_3-x_1 \rangle}{\R^{3d}}{\Bmu} + o(W_{2,\Bmu}(\mu^1,\mu^3)),
\end{equation}
for all $\mu^3 \in \Pcal_2(\R^d)$ and $\Bmu \in \Gamma(\Bgamma,\mu^3)$. Here for $\Bmu \in \Pcal(\R^{3d})$, the quantity $W_{2,\Bmu}(\mu^1,\mu^3)$ is defined by
\begin{equation*}
W_{2,\Bmu}(\mu^1,\mu^3) = \left( \INTDom{|x_1-x_3|^2}{\R^{2d}}{\Bmu(x_1,x_2,x_3)} \right)^{1/2}.
\end{equation*}
\end{Def}

We now introduce the technical notions of \textit{regularity} and \textit{metric slope} that are instrumental in deriving a sufficient condition for the extended subdifferential of a functional to be non-empty. This result is stated in Theorem \ref{thm:subdifferential_localslope} and its proof can be found in \cite[Theorem 10.3.10]{AGS}.

\begin{Def}[Regular functionals over $(\Pcal_2(\R^d),W_2)$ and metric slope]
\label{def:Regular}
A proper and lower semicontinuous functional $\phi(\cdot)$ is said to be \textit{regular} provided that whenever $(\mu_n) \subset \Pcal_2(\R^d)$ and $(\Bgamma_n) \subset \Pcal_2(\R^{2d})$ are taken such that
\begin{equation*}
\left\{
\begin{aligned}
& \mu_n \overset{W_2 \,}{\longrightarrow} \mu ~~ \text{in $\Pcal_2(\R^d)$} ~~,~~ \phi(\mu_n) \longrightarrow \tilde{\phi} ~~ \text{in $\R$}, \\
& \Bgamma_n \in \Bpartial_S \phi(\mu_n) ~ \forall n \geq 1 ~,~ \Bgamma_n \overset{W_2 \,}{\longrightarrow} \Bgamma ~ \text{in $\Pcal_2(\R^{2d})$}, 
\end{aligned}
\right.
\end{equation*}
it implies that $\Bgamma \in \Bpartial \phi(\mu)$ and $\tilde{\phi} = \phi(\mu)$.

Furthermore, we define the \textit{metric slope} $|\partial \phi|(\mu)$ of the functional $\phi(\cdot)$ at $\mu \in D(\phi)$ as
\begin{equation*}
| \partial \phi |(\mu) = \underset{\nu \rightarrow \mu}{\textnormal{limsup}} \left[ \frac{\left( \phi(\mu) - \phi(\nu) \right)^+}{W_2(\mu,\nu)}\right].
\end{equation*}
where $(\bullet)^+$ denotes the positive part.
\end{Def}

\begin{thm}[Link between extended subdifferentials and metric slopes]
\label{thm:subdifferential_localslope}
Let $\phi(\cdot)$ be a proper, lower-semicontinuous, bounded from below and regular functional over $\Pcal_2(\R^d)$. Then, the extended subdifferential $\Bpartial \phi(\mu)$ of $\phi(\cdot)$ at some $\mu \in D(\phi)$ is non-empty if  and only if its metric slope $|\partial \phi|(\mu)$ at $\mu$ is finite.

In which case, there exists a unique \textnormal{minimal selection} in $\Bpartial \phi(\mu)$, denoted by $\Bpartial^{\circ} \phi(\mu)$, satisfying
\begin{equation*}
\begin{aligned}
\left( \INTDom{|r|^2}{\R^{2d}}{(\Bpartial^{\circ} \phi(\mu))(x,r)} \right)^{1/2} = & \min \left\{ \left( \INTDom{|r|^2}{\R^{2d}}{\Bgamma(x,r)} \right)^{1/2} ~ \text{s.t.} ~ \Bgamma \in \Bpartial \phi(\mu) \right\} = |\partial \phi |(\mu).
\end{aligned}
\end{equation*}
The minimal selection can be explicitly characterized as follows: let $\mu_{\tau}$ be the minimizer of the Moreau-Yosida functional \eqref{eq:Moreau_Yosida} for some $\tau \in (0,\tau_*)$. Then there exists a family of strong subdifferentials $(\Bgamma_{\tau}) \subset (\Bpartial_S \phi(\mu_{\tau}))$ which converges towards $\Bpartial^{\circ} \phi(\mu)$ in the $W_2$-metric along any vanishing sequence $\tau_n \downarrow 0$.
\end{thm}

We list in Section \ref{section:Examples} below several examples of regular functionals and compute the minimal selection in their extended subdifferential. We end this section by recalling the definition of \textit{Wasserstein gradient}.

\begin{Def}[Wasserstein gradient]
\label{def:Wass_Grad}
Whenever the minimal selection $\Bpartial^{\circ} \phi(\mu)$ is induced by a Borel map, this map is called the \textit{Wasserstein gradient} of $\phi(\cdot)$. It is denoted by $\nabla_{\mu} \phi(\mu) \in L^2(\R^d,\R^d;\mu)$ and it coincides with the barycenter of the minimal selection.
\end{Def}

The main interest of subdifferential calculus in the space $(\Pcal_2(\R^d),W_2)$ is to compute derivatives of functionals along measure curves. However, the general chain rule described in \cite[Proposition 10.3.18]{AGS} only applies to the case of a curve $\epsilon \mapsto \mu(\epsilon) = \G(\epsilon,\cdot)_{\#}\mu$ generated by a given smooth functions $\G(\epsilon,\cdot)$ when one restricts himself to strong subdifferentials. Yet, there is no reason in general for the strong subdifferential of a functional to be non-empty. 
In Proposition \ref{prop:Chainrule}, we condense some well known results of \cite[Chapter 10]{AGS} in order to provide a chain rule that allows to compute derivatives along smooth vector fields using the minimal selection $\Bpartial^{\circ} \phi(\mu)$. For simplicity, we state this result in the framework of the Wasserstein space $\Pcal_c(\R^d)$.

\begin{prop}[Minimal selection and chain rule along smooth vector fields]
\label{prop:Chainrule}
Let $\mu \in \Pcal_c(\R^d)$ and $K = \cup_{x \in \supp(\mu)} \overline{B(x,1)}$. Let $\phi : \Pcal(K) \rightarrow (-\infty,+\infty]$ be a functional satisfying hypotheses \textnormal{\textbf{(C)}} and \textnormal{\textbf{(D)}} of Theorem \ref{thm:General_PMP}. Define $\G \in \Lip((-\bar{\epsilon},\bar{\epsilon}),C^0(\R^d,\R^d))$ a family of continuous functions with $\G(0,\cdot) = \Id$, $\supp(\G(\epsilon,\cdot)_{\#}\mu) \subset K$ for all $\epsilon \in (-\bar{\epsilon},\bar{\epsilon})$ and $\F : x \mapsto \derv{}{\epsilon}{} \left[ \G(\epsilon,x) \right]_{\epsilon = 0}$ being $C^0$ as well.

Then it holds that 
\begin{equation*}
\derv{}{\epsilon}{} \left[ \phi(\G(\epsilon,\cdot)_{\#}\mu) \right]_{\epsilon = 0} = \INTDom{\langle \bar{\Bgamma}^{\circ}(x) , \F(x) \rangle}{\R^{2d}}{\mu(x)},
\end{equation*}
where $\bar{\Bgamma}^{\circ} \in L^2(\R^d,\R^d;\mu)$ is the barycenter of $\Bpartial^{\circ} \phi(\mu)$.
\end{prop}

\begin{proof}
First remark that it holds for any $\nu \in \Pcal(K)$
\begin{equation*}
\left( \phi(\mu) - \phi(\nu) \right)^+ \leq \Lip(\phi,\Pcal(K)) W_2(\mu,\nu)\end{equation*}
where $\Lip(\phi,\Pcal(K))$ is the Lipschitz constant of $\phi(\cdot)$ on $\Pcal(K)$. Hence, $|\partial \phi|(\mu)$ is uniformly bounded by $\Lip(\phi,\Pcal(K))$. Moreover, the assumption that $\phi(\cdot)$ is bounded from below and Lipschitz on sets of uniformly compactly supported measures implies that for $\tau_* > 0$ small enough, the Moreau-Yosida functional $\Phi_{\text{M}}(\mu,\tau;\cdot)$ defined in \eqref{eq:Moreau_Yosida} attains a minimum point $\mu_{\tau} \in D(\phi) \subset \Pcal(K)$ for any $\tau \in (0,\tau_*)$. Thus, by Theorem \ref{thm:subdifferential_localslope}, $\Bpartial \phi(\mu)$ is non-empty and contains at least the minimal selection $\Bpartial^{\circ} \phi(\mu)$ at any $\mu \in \Pcal(K)$.

Consider a sequence $(\tau_n) \subset (0,\tau_*)$ converging to $0$ and the corresponding sequence of strong subdifferentials $(\Bgamma_{\tau_n}) \subset (\Bpartial_S \phi(\mu_{\tau_n}))$ converging towards $\Bpartial^{\circ} \phi(\mu)$ in the $W_2$-metric. Pick $\epsilon \in (0,\bar{\epsilon})$ small enough and choose $\Bmu_{\epsilon}^{\tau_n} = (\pi^1 , \pi^2 , \G(\epsilon,\cdot) \circ \pi^1)_{\#} \Bgamma_{\tau_n} \in \Gamma(\Bgamma_{\tau_n},\G(\epsilon,\cdot)_{\#}\mu_{\tau_n})$. By the definition of strong subdifferentials given in \eqref{eq:Def_StrongSubdiff}, it holds that
\begin{equation} \label{eq:Proof_Chainrule_Ineq}
\frac{\phi(\G(\epsilon,\cdot)_{\#}\mu_{\tau_n}) - \phi(\mu_{\tau_n})}{\epsilon} \geq \INTDom{\langle r , \frac{\G(\epsilon,x) - x}{\epsilon} \rangle}{\R^{2d}}{\Bgamma_{\tau_n}(x,r)} + o(1).
\end{equation}
since 
\begin{equation*}
o(W_{2,\Bmu}(\G(\epsilon,\cdot)_{\#} \mu_{\tau_n} , \mu_{\tau_n})) = o \left( \NormL{\G(\epsilon,\cdot)-\Id}{2}{\mu_{\tau_n}} \right) = o(\epsilon) ~ \text{for all $n \geq 1$}.
\end{equation*}
Remark that the left hand side of \eqref{eq:Proof_Chainrule_Ineq} is bounded over $\Pcal(K)$ uniformly with respect to $n \geq 1$ and $\epsilon \in (0,\bar{\epsilon})$ by Lipschitzianity of $\phi(\cdot)$.

We recall that $\Bgamma_{\tau_n} \overset{W_2 \,}{\longrightarrow} \Bpartial^{\circ} \phi(\mu)$ in $\Pcal_2(\R^{2d})$. Notice that the whole sequence $(\mu_{\tau_n})$ is in $\Pcal(K)$, thus for all $\epsilon \in (0,\bar{\epsilon})$ the maps $x \mapsto |(\G(\epsilon,x)-x)/\epsilon|^2$ are uniformly integrable with respect to $\{ \pi^1_{\#} \Bgamma_{\tau_n} \}_{n=1}^{+\infty}$. Hence, the maps $(x,r) \mapsto \left| \langle r , (\G(\epsilon,x)-x)/\epsilon \rangle \right|$ are uniformly integrable with respect to $\{ \Bgamma_n \}_{n=1}^{+\infty}$ and the application of Proposition \ref{prop:Convergence_results}-$(i)$ implies that for all $\epsilon \in (0,\bar{\epsilon})$,
\begin{equation}
\label{eq:tau_nConvergence}
\begin{aligned}
\INTDom{\langle r , \frac{\G(\epsilon,x) - x}{\epsilon} \rangle}{\R^{2d}}{\Bgamma_{\tau_n}(x,r)} ~ \underset{\tau_n \downarrow 0}{\longrightarrow} ~ & \INTDom{\langle r , \frac{\G(\epsilon,x) - x}{\epsilon} \rangle}{\R^{2d}}{(\Bpartial^{\circ}\phi(\mu))(x,r)} 
\\ = \hspace{0.4cm} & \INTDom{\langle \bar{\Bgamma}^{\circ}(x) , \frac{\G(\epsilon,x) - x}{\epsilon} \rangle}{\R^d}{\mu(x)}
\end{aligned}
\end{equation}
using the notion of \textit{barycenter} of a plan introduced in Definition \ref{def:Barycenter}. 

Moreover, the Lipschitz regularity in the $W_2$-metric of $\phi(\cdot)$ over $\Pcal(K)$ together with Proposition \ref{prop:Convergence_results}-$(ii)$ imply that
\begin{equation}
\label{eq:phiConvergence}
\phi(\G(\epsilon,\cdot)_{\#} \mu_{\tau_n}) \longrightarrow \phi(\G(\epsilon,\cdot)_{\#} \mu).
\end{equation}
Thus, merging \eqref{eq:Proof_Chainrule_Ineq},\eqref{eq:tau_nConvergence} and \eqref{eq:phiConvergence}, we prove that for any $\epsilon \in (0,\bar{\epsilon})$ with $\bar{\epsilon} > 0$ small enough, it holds 
\begin{equation*}
\frac{\phi(\G(\epsilon,\cdot)_{\#} \mu) - \phi(\mu)}{\epsilon} \geq \INTDom{\langle \bar{\Bgamma}^{\circ}(x) , \frac{\G(\epsilon,x) - x}{\epsilon} \rangle}{\R^d}{\mu(x)} + o(1).
\end{equation*}

Invoking similar arguments, the family of maps $( \left| \langle \bar{\Bgamma}^{\circ}(\cdot) , (\G(\epsilon,\cdot) -\Id)/\epsilon \rangle \right| )_{\epsilon \in (0,\bar{\epsilon})}$ is uniformly integrable with respect to $\mu$ and it holds that $\left| \langle \bar{\Bgamma}^{\circ}(\cdot) , \F(\cdot) \rangle \right| < +\infty$ $\mu$-almost everywhere. Therefore, letting $\epsilon \downarrow 0$ and invoking Proposition \ref{prop:Convergence_results}-$(iii)$, we recover that 
\begin{equation*}
\lim\limits_{\epsilon \downarrow 0} \left[ \frac{\phi(\G(\epsilon,\cdot)_{\#} \mu) - \phi(\mu)}{\epsilon} \right] \geq \INTDom{\langle \bar{\Bgamma}^{\circ}(x) , \F(x) \rangle}{\R^d}{\mu(x)}.
\end{equation*}

Following the same steps with $\epsilon \in (-\bar{\epsilon},0)$, we obtain the converse inequality for $\epsilon \uparrow 0$. Since we assumed that $\epsilon \mapsto \phi(\G(\epsilon,\cdot)_{\#} \mu)$ is differentiable at $\epsilon = 0$ in \textbf{(D)}, these limits coincide and it holds 
\begin{equation*}
\derv{}{\epsilon}{}\left[ \phi(\G(\epsilon,\cdot)_{\#}\mu) \right]_{\epsilon = 0} = \INTDom{\langle \bar{\Bgamma}^{\circ}(x) , \F(x) \rangle}{\R^d}{\mu(x)},
\end{equation*} 
which proves our claim. 
\end{proof}

\begin{rmk}[The case $\Bpartial_S \phi(\mu) \neq \emptyset$]
When $\Bpartial_S \phi(\mu)$ is non-empty, the previous chain rule can be applied with any strong subdifferential and for more general classes of vector fields, see e.g. \cite[Remark 10.3.2]{AGS}.
\end{rmk}

The interest of proving this kind of result for the minimal selection is twofold. First, as recalled in Theorem \ref{thm:subdifferential_localslope}, a minimal selection always exists when the extended subdifferential is non-empty. Second, minimal selections can be computed explicitly even in very general settings for a wide range of functionals (see e.g. \cite[Chapter 10.4]{AGS} or Section \ref{section:Examples}). In such cases, they are usually induced by their barycenter, yielding the existence of a \textnormal{Wasserstein gradient} for the functional.

%%%%%%%%%%%%%%%%%%%%%%%%%%%%%%%%%%%%%%%%%%%%%%%%%%%%%%%%%%%%%%%%%%%%%%%%%%%%%%%%%%

\subsection{Directional derivatives of non-local flows}

In this section, we prove the existence of directional derivatives along measure curves generated by suitable Lipschitz families of continuous and bounded maps for non-local flows. Such derivatives are characterized as the only solution of a linearized Cauchy problem. This result can be seen as a generalization to the Wasserstein setting of Proposition \ref{prop:Differential_local_flow}.

Before stating our result, we recall the classical Banach Fixed Point Theorem with parameter (see e.g. \cite[Theorem A.2.1]{BressanPiccoli}).

\begin{thm}[Banach fixed point theorem with parameter]
\label{thm:Banach}
Let $X$ be a Banach space, $S$ be a metric space and $\Lambda : X \times S \rightarrow X$ be a continuous mapping such that, for some $\kappa < 1$,
\begin{equation*}
\| \Lambda(x,s) - \Lambda(y,s) \|_X \leq \kappa \| x-y \|_X ~~ \text{for all $x,y \in X$ and $s \in S$.}
\end{equation*}
Then for each $s \in S$, there exists a unique fixed point $x(s) \in X$ of $\Lambda(\cdot,s)$. Moreover, the map $s \mapsto x(s)$ is continuous and for any $(s,y) \in S \times X$, it holds
\begin{equation}
\label{eq:Parametrized_FixedPoint}
\| y - x(s) \|_X \leq \frac{1}{1-\kappa} \|y - \Lambda(s,y) \|_X.
\end{equation}
\end{thm}

We are now ready to state and prove the main result of this Section.

\begin{prop}[Directional derivative of a non-local flow with respect to the initial data]
\label{prop:Directional_derivative_flow} 
Let $\mu \in \Pcal_c(\R^d)$, $\bar{\epsilon} > 0$ be a small parameter, $\G(\cdot,\cdot) \in \Lip((-\bar{\epsilon},\bar{\epsilon}),C^0(\R^d,\R^d))$ be a family of bounded maps with $\G(0,\cdot) = \Id$ and $\F : x \in \supp(\mu) \mapsto \derv{}{\epsilon}{} \left[ \G(\epsilon,x) \right]_{\epsilon=0}$ be continuous as well. 

Let $v[\cdot](\cdot,\cdot)$ be a non-local vector field satisfying hypotheses \textnormal{\textbf{(F),(B),(D)}}, $\Phi^v_{(0,\cdot)}[\cdot](\cdot)$ be the corresponding family of non-local flows as defined in Theorem \ref{thm:Existence_uniqueness_non-local_PDE} and $\mu(\cdot)$ be the unique solution of the corresponding Cauchy problem \eqref{eq:Nonlocal_Cauchy_problem} starting from $\mu$.

Then, the map $\epsilon \in (-\bar{\epsilon},\bar{\epsilon}) \mapsto \Phi^v_{(0,t)}[\G(\epsilon,\cdot)_{\#}\mu](x)$ admits a derivative at $\epsilon = 0$ for all $(t,x) \in [0,T] \times \overline{B(0,R_T)}$ that we denote by $w_{\Phi}(t,x)$. It can be characterised as the unique solution of the Cauchy problem
\begin{equation}
\label{eq:Directional_derivative_charac}
\left\{
\begin{aligned}
& \partial_t w(t,x) = \D_x v[\mu(t)] \left( t, \Phi^v_{(0,t)}[\mu](x) \right) w(t,x) \\
& + \INTDom{ \BGamma^{\circ}_{\left( t, \Phi^v_{(0,t)}[\mu](x) \right)} \left( \Phi^v_{(0,t)}[\mu](y) \right) \cdot \left[ \D_x \Phi_{(0,t)}[\mu](y) \F(y) + w(t,y) \right]  }{\R^d}{\mu(y)}, \\
& w(0,x) = 0 \hspace{7.375cm} \text{for all $x \in \R^d$},
\end{aligned}
\right.
\end{equation}
where for all $(t,z)$, $\BGamma^{\circ}_{\left( t,z \right)}(\cdot)$ is the matrix-valued map made of the barycenters of the minimal selections $\Bpartial_{\mu}^{\circ} v^i[\mu(t)](t,z)$ in the extended subdifferential of the components of $\mu \mapsto v^i[\mu](t,z)$ at $\mu(t)$.
\end{prop}

\begin{proof}
We follow a classical scheme of proof used in the finite dimensional setting to show that flows of diffeomorphims admit directional derivatives characterized as the unique solution of a linearized Cauchy problem (see e.g. \cite[Theorem 2.3.1]{BressanPiccoli}. 

First, we define $\Omega = \overline{B(0,R_T)}$ and we introduce the operator $\Lambda_{\Phi} : w \in C^0([0,T] \times \Omega,\R^d) \mapsto \Lambda_{\Phi}(w) \in C^0([0,T] \times \Omega,\R^d)$ defined for all $(t,x) \in [0,T] \times \Omega$ by 
\begin{equation*}
\begin{aligned}
& \Lambda_{\Phi}(w) (t,x) = \INTSeg{ \D_x v[\mu(s)] \left( s , \Phi^v_{(0,s)}[\mu](x) \right) w(s,x)}{s}{0}{t} \\
 + & \INTSeg{\INTDom{ \BGamma^{\circ}_{\left( s , \Phi^v_{(0,s)}[\mu](x) \right)} \left( \Phi^v_{(0,s)}[\mu](y) \right) \cdot \left[ \D_x \Phi_{(0,s)}[\mu](y) \F(y) + w(s,y) \right]  }{\R^d}{\mu(y)}}{s}{0}{t}.
\end{aligned}
\end{equation*}
By hypotheses \textbf{(F)} and \textbf{(B)}, the right hand side of the previous equation is continuous in $(t,x)$. We first show that this operator admits a unique fixed point and afterwards that it coincides with the map which to every $(t,x)$ associates the derivative at $\epsilon = 0$ of the family of non-local flows $\epsilon \mapsto \Phi^v_{(0,t)}[\G(\epsilon,\cdot)_{\#}\mu](x)$. With this goal, we introduce a parameter $\alpha > 0$ that will be chosen so that the operator $\Lambda_{\Phi}(\cdot)$ is contracting with respect to the equivalent norm
\begin{equation}
\label{eq:Equivalent_C0norm}
\NormC{w}{0}{[0,T] \times \Omega}^{\alpha} = \sup\limits_{(t,x) \in [0,T] \times \Omega} \left[ e^{-2\alpha t} |w(t,x)| \right].
\end{equation}

Remark that for any $w_1,w_2 \in C^0([0,T] \times \Omega,\R^d)$ and any $(t,x) \in [0,T] \times \Omega$, it holds
\begin{equation*}
\begin{aligned}
| \Lambda_{\Phi}(w_2)(t,x) - \Lambda_{\Phi}(w_1)(t,x) | & \leq \INTSeg{ \left| \D_x v[\mu(s)] \left( s , \Phi^v_{(0,s)}[\mu](x) \right) \cdot (w_2(s,x)-w_1(s,x)) \right| }{s}{0}{t} \\
& + \INTSeg{ \INTDom{ \left| \BGamma^{\circ}_{\left( s , \Phi^v_{(0,s)}[\mu](x) \right)} \left( \Phi^v_{(0,s)}[\mu](y) \right) \cdot ( w_2(s,y) - w_1(s,y) ) \right| }{\R^d}{\mu(y)}}{s}{0}{t} \\
& \leq \INTSeg{ \left( L_1 |w_2(s,x) - w_1(s,x)| + L_2 \NormL{w_2(s,\cdot)-w_1(s,\cdot)}{1}{\mu} \right)}{s}{0}{t} \\
& \leq \INTSeg{ (L_1 + L_2) \NormC{w_2(s,\cdot)-w_1(s,\cdot)}{0}{\Omega}}{s}{0}{t},
\end{aligned}
\end{equation*}
since $\mu(\Omega) = 1$, and where we introduced
\begin{equation*}
L_1 = \left\| \D_x v[\mu(\cdot)] (\cdot , \Phi^v_{(0,\cdot)}[\mu](\cdot) ) \right\|_{L^{\infty}([0,T] \times \Omega ; \Lcal^1 \times \mu(\cdot))},
\end{equation*}
and
\begin{equation*}
L_2 = \left\| \BGamma^{\circ}_{\left( \cdot , \Phi^v_{(0,\cdot)}[\mu](\cdot) \right)} \left( \Phi^v_{(0,\cdot)}[\mu](\cdot) \right)\right\|_{L^{\infty}([0,T] \times \Omega^2 ; \Lcal^1 \times \mu(\cdot) \times \mu(\cdot) )}
\end{equation*}
which exist by hypotheses \textbf{(F)} and \textbf{(B)}. It further holds by definition of $\NormC{\cdot}{0}{[0,T]\times \Omega}^{\alpha}$ that
\begin{equation*}
\begin{aligned}
| \Lambda_{\Phi}(w_2)(t,x) - \Lambda_{\Phi}(w_1)(t,x) | & \leq \INTSeg{ e^{2 \alpha s} (L_1 + L_2) \NormC{w_2(\cdot,\cdot) -w_1(\cdot,\cdot)}{0}{[0,T] \times \Omega}^{\alpha}}{s}{0}{t} \\
& \leq \frac{e^{2 \alpha t}-1}{2 \alpha}(L_1 + L_2) \NormC{w_2(\cdot,\cdot) -w_1(\cdot,\cdot)}{0}{[0,T] \times \Omega}^{\alpha}.
\end{aligned}
\end{equation*}
Multiplying both sides of the inequality by $e^{-2\alpha t}$ and taking the supremum over $(t,x) \in [0,T] \times \Omega$ in the left-hand side yields the desired contractivity with a constant equal to $1/2$ provided that $\alpha \geq (L_1 + L_2)$. It is then possible to apply Theorem \ref{thm:Banach} to obtain the existence of a unique fixed point $w_{\Phi}(\cdot,\cdot) \in C^0([0,T]\times \Omega , \R^d)$ of $\Lambda_{\Phi}(\cdot)$.  

Define for $\epsilon \in (-\bar{\epsilon},\bar{\epsilon})$ the parametrized family of operators $\Psi^{\epsilon} : f \in  C^0([0,T] \times \Omega , \R^d)  \mapsto  \Psi^{\epsilon}(f) \in C^0([0,T] \times \Omega , \R^d)$ defined by
\begin{equation}
\label{eq:Def_psi_epsilon}
\Psi^{\epsilon}(f)(t,x) = x + \INTSeg{v[f(s,\cdot)_{\#} (\G(\epsilon,\cdot)_{\#} \mu)] (s,f(s,x))  }{s}{0}{t}
\end{equation}
for all $(t,x) \in [0,T] \times \Omega$. Up to defining again the equivalent norm $\NormC{\cdot}{0}{[0,T] \times \Omega}^{\alpha}$ as in \eqref{eq:Equivalent_C0norm} with a suitable $\alpha > 0$, it can be shown that this operator is contracting independently from $\epsilon$ as a direct consequence of the Lipschitzianity hypotheses given in \textbf{(F)}. We can thus invoke again Theorem \ref{thm:Banach} to obtain the existence of a unique fixed point of $\Psi^{\epsilon}(\cdot)$ for each $\epsilon \in (-\bar{\epsilon},\bar{\epsilon})$. Notice that by definition, this family of fixed points is precisely the parametrized family of non-local flows $(t,x) \mapsto \Phi^v_{(0,t)}[\G(\epsilon,\cdot)_{\#}\mu](x)$. As a consequence of Theorem \ref{thm:Existence_uniqueness_non-local_PDE}, we know that these maps are $C^1$ with respect to $x$ for all $\epsilon$.

We now define the map $\hat{\Phi}^{v,\epsilon}_{(0,\cdot)}[\mu](\cdot)$ by 
\begin{equation*}
\hat{\Phi}^{v,\epsilon}_{(0,\cdot)}[\mu](\cdot) : (t,x) \mapsto \Phi^v_{(0,t)}[\mu](x) + \epsilon w_{\Phi}(t,x).
\end{equation*}
To conclude, we then need to show that
\begin{equation*}
\lim\limits_{\epsilon \rightarrow 0} \NormC{ \frac{1}{\epsilon} \left(  \hat{\Phi}^{v,\epsilon}_{(0,\cdot)}[\mu](\cdot) - \Phi^v_{(0,\cdot)}[\G(\epsilon,\cdot)_{\#}\mu](\cdot) \right) }{0}{[0,T] \times \Omega}  = 0,
\end{equation*}
which will directly yield the existence and the characterization of the directional derivative of the flow along $(-\bar{\epsilon},\bar{\epsilon}) \mapsto \G(\epsilon,\cdot)_{\#} \mu^0$. By \eqref{eq:Parametrized_FixedPoint} in Theorem \ref{thm:Banach} and the equivalence of the $C^0$-norms we introduced, there exists a constant $C > 0$ independent from $\epsilon$ such that it holds
\begin{equation*}
\frac{1}{|\epsilon|}\left \|  \hat{\Phi}^{v,\epsilon}_{(0,\cdot)}[\mu](\cdot) -\Phi^v_{(0,\cdot)}[\G(\epsilon,\cdot)_{\#}\mu](\cdot) \right \|_{C^0} \leq \frac{2C}{| \epsilon |} \left\| \hat{\Phi}^{v,\epsilon}_{(0,\cdot)}[\mu](\cdot) -\Psi^{\epsilon}(\hat{\Phi}^{v,\epsilon}_{(0,\cdot)}[\mu](\cdot))(\cdot,\cdot)\right \|_{C^0}.
\end{equation*}
We now want to perform a first order expansion on $\Psi^{\epsilon}(\hat{\Phi}^{v,\epsilon}_{(0,s)}[\mu](\cdot))(\cdot,\cdot)$ with respect to $\epsilon$. Take $(s,x) \in [0,T] \times \Omega$. One has by definition of $\hat{\Phi}^{v,\epsilon}_{(0,\cdot)}[\mu](\cdot)$ that 
\begin{equation*}
\begin{aligned}
\hat{\Phi}^{v,\epsilon}_{(0,s)}[\mu](\G(\epsilon,x)) = & \, \Phi^v_{(0,s)}[\mu](\G(\epsilon,x)) + \epsilon w_{\Phi}(s,\G(\epsilon,x)) \\
= & \, \Phi^v_{(0,s)}[\mu](x) + \epsilon \left( \D_x \Phi^v_{(0,s)}[\mu](x) \F(x) + w_{\Phi}(s,x) \right) + o(\epsilon),
\end{aligned}
\end{equation*}
by continuity of $w_{\Phi}(s,\cdot)$ for all $s \in [0,T]$.

By assumptions \textbf{(F)}, \textbf{(B)} and \textbf{(D)}, we can apply the chain rule of Proposition \ref{prop:Chainrule} component-wise on the $v^i$ to obtain that
\begin{equation} 
\label{eq:DirectionalDerivative_Est1}
\begin{aligned}
& v \left[ \hat{\Phi}^{v,\epsilon}_{(0,s)}[\mu](\cdot) \circ \G(\epsilon,\cdot)_{\#} \mu \right](s,z) = v[\Phi^v_{(0,s)}[\mu](\cdot)_{\#}\mu](s,z) \\
+ \, & \epsilon \INTDom{\BGamma^{\circ}_{(s,z)} \left( \Phi^v_{(0,s)}[\mu](y) \right) \left[ \D_x \Phi^v_{(0,s)}[\mu](y) \F(y) + w_{\Phi}(s,y) \right] }{\R^d}{\mu(y)} + o(\epsilon)
\end{aligned}
\end{equation}
where for all $(s,z)$ the map $y \mapsto \BGamma^{\circ}_{(s,z)}(y) = (\bar{\Bgamma}^{i,\circ}_{(s,z)}(y))_{1 \leq i \leq d} \in \R^{d \times d}$ is made of the barycenters of the minimal selections in the extended subdifferentials of the components $v^i$'s. 

Performing a Taylor expansion in the space variable for the non-local velocity field, it also holds that
\begin{equation}
\label{eq:DirectionalDerivative_Est2}
\begin{aligned}
v \left[ \Phi^v_{(0,s)}[\mu](\cdot)_{\#}\mu \right] \left(s,\Phi^v_{(0,s)}[\mu](x) + \epsilon w_{\Phi}(s,x) \right) & = v[\Phi^v_{(0,s)}[\mu](\cdot)_{\#} \mu] \left( s,\Phi^v_{(0,s)}[\mu](x) \right) \\ 
& + \epsilon \D_x v[\Phi^v_{(0,s)}[\mu](\cdot)_{\#}\mu] \left( s,\Phi^v_{(0,s)}[\mu](x) \right) \cdot w_{\Phi}(s,x) + o(\epsilon),
\end{aligned}
\end{equation}
as well as 
\begin{equation}
\label{eq:DirectionalDerivative_Est3}
\BGamma^{\circ}_{\left( s, \Phi^v_{(0,s)}[\mu](x) + \epsilon w_{\Phi}(s,x) \right)} \left( \Phi^v_{(0,s)}[\mu](y) \right) = \BGamma^{\circ}_{\left( s, \Phi^v_{(0,s)}[\mu](x) \right)} \left( \Phi^v_{(0,s)}[\mu](y) \right) + o(1) 
\end{equation}
thanks to assumption \textbf{(B)} in which we state that $z \mapsto \BGamma^{\circ}_{(t,z)}(x)$ is continuous for all $(s,x) \in [0,T] \times \R^d$ .

Merging \eqref{eq:Def_psi_epsilon}, \eqref{eq:DirectionalDerivative_Est1} and \eqref{eq:DirectionalDerivative_Est2}, \eqref{eq:DirectionalDerivative_Est3} and recalling the definition of $w_{\Phi}(\cdot,\cdot)$, it holds 
\begin{equation*}
\begin{aligned}
& \Psi^{\epsilon} \left( \hat{\Phi}^{v,\epsilon}_{(0,\cdot)}[\mu](\cdot) \right)(t,x) = x + \INTSeg{v[\Phi^v_{(0,s)}[\mu](\cdot)_{\#}\mu](s,\Phi^v_{(0,s)}[\mu](x))}{s}{0}{t} \\ 
& \hspace{0.45cm} + \epsilon \INTSeg{\INTDom{\BGamma^{\circ}_{\left( s,\Phi^v_{(0,s)}[\mu](x) \right)} \left( \Phi^v_{(0,s)}[\mu](y) \right) \left[ \D_x \Phi^v_{(0,s)}[\mu](y) \F(y) + w_{\Phi}(s,y) \right] }{\R^d}{\mu(y)}}{s}{0}{t} \\
& \hspace{0.45cm} + \epsilon \INTSeg{\D_x v[\Phi^v_{(0,s)}[\mu](\cdot)_{\#}\mu](s,\Phi^v_{(0,s)}[\mu](x)) \cdot w_{\Phi}(s,x)}{s}{0}{t} + o(\epsilon) \\
=  ~ & \Phi^{v}_{(0,t)}[\mu](x) + \epsilon w_{\Phi}(t,x) + o(\epsilon).
\end{aligned}
\end{equation*}

Therefore, we finally recover that
\begin{equation*}
\frac{1}{| \epsilon |} \left| \Psi^{\epsilon}(\hat{\Phi}^{\epsilon}_{(0,\cdot)}[\mu](\cdot))(t,x) - \hat{\Phi}^{\epsilon}_{(0,t)}[\mu](x) \right| \leq o(1)
\end{equation*}
as $\epsilon \rightarrow 0$ for all $(t,x) \in [0,T] \times \Omega$, and conclude that
\begin{equation*}
\lim\limits_{\epsilon \rightarrow 0} \left[ \NormC{ \frac{1}{\epsilon} \left( \Phi_{(0,\cdot)}[\G(\epsilon,\cdot)_{\#}\mu](\cdot) - \hat{\Phi}^{\epsilon}_{(0,\cdot)}[\mu](\cdot) \right) }{0}{[0,T] \times \Omega} \right] = 0.
\end{equation*}

We thus proved that the derivative of $\epsilon \in (-\bar{\epsilon},\bar{\epsilon}) \mapsto \Phi^v_{(0,t)}[\G(\epsilon,\cdot)_{\#}\mu](x)$ at $\epsilon = 0$ exists for any $(t,x)$ and that it is the only solution of equation \eqref{eq:Directional_derivative_charac}.
\end{proof}

%%%%%%%%%%%%%%%%%%%%%%%%%%%%%%%%%%%%%%%%%%%%%%%%%%%%%%%%%%%%%%%%%%%%%%%%%%%%
%																		   
%									NEW SECTION
%
%%%%%%%%%%%%%%%%%%%%%%%%%%%%%%%%%%%%%%%%%%%%%%%%%%%%%%%%%%%%%%%%%%%%%%%%%%%%

\section{The Pontryagin Maximum Principle}
\label{section:PMP}

In this section, we state the main result of our article.

\begin{thm}[Pontryagin Maximum Principle for $(\pazocal{P})$] 
\label{thm:General_PMP}
Let $(u^*(\cdot),\mu^*(\cdot)) \in \U \times \Lip([0,T],\Pcal_c(\R^d))$ be an optimal pair control-trajectory for $(\pazocal{P})$ and assume that the following hypotheses \textnormal{\textbf{(H)}} hold :

\begin{framed} 
\vspace{-0.2cm}
\begin{center} \textnormal{\textbf{(H)}} \end{center}
\begin{enumerate}
\item[] \textnormal{\textbf{(U)}} : The set of admissible controls is $\U = L^1([0,T],U)$ where $U \subset C^1(\R^d,\R^d)$ is a non-empty and closed subset of $\{ v \in C^1(\R^d,\R^d) ~\text{s.t.}~ \NormC{v(\cdot)}{1}{\R^d} \leq L_U \}$ for a given constant $L_U > 0$. \smallskip
\item[] \textnormal{\textbf{(L)}} : The running cost $L : (\mu,\omega) \in \Pcal_c(\R^d) \times U \mapsto L(\mu,\omega) \in \R$ is Lipschitz in $(\mu,\omega)$ with respect to the product metric $W_2 \times C^0$ over $\Pcal(K) \times U$ for any compact set $K \subset \R^d$. The functional $\mu \in \Pcal(K) \mapsto L(\mu,\omega)$ is proper, regular in the sense of Definition \ref{def:Regular} below, bounded for any $\omega \in U$ and $K \subset \R^d$ compact. \smallskip
\item[] \textnormal{\textbf{(C)}} : The terminal cost $\varphi : \mu \in \Pcal_c(\R^d) \mapsto \varphi(\mu) \in \R$ is proper, regular in the sense of Definition \ref{def:Regular} below, Lipschitz with respect to the $W_2$-metric, bounded from below over $\Pcal(K)$ for any compact set $K \subset \R^d$.  \smallskip
\item[] \textnormal{\textbf{(F)}} : The non-local velocity field $v : \mu \in \Pcal_c(\R^d) \mapsto v[\mu](\cdot,\cdot) \in L^1([0,T],C^1(\R^d,\R^d) \cap L^{\infty}(\R^d,\R^d))$ satisfies
\begin{equation*}
\begin{aligned}
& | v[\mu](t,x) | \leq M(1+|x|) ~,~ | v[\mu](t,x) - v[\mu](t,y) | \leq L_1|x-y| ~, \\
\text{and} ~~ & \hspace{-0.1cm} \NormC{v[\mu](t,\cdot) - v[\nu](t,\cdot)}{0}{\R^d} \leq L_2 W_1(\mu,\nu) 
\end{aligned}
\end{equation*}
for $\Lcal^1$-almost every $t \in [0,T]$ and all $(x,y) \in  \R^{2d}$ where $M,L_1$ and $L_2$ are positive constants. For any compact set $K \subset \R^d$ and any $i \in \{ 1,...,d \}$, the components $\mu \in \Pcal(K) \mapsto v^i[\mu](t,x)$ are regular in the sense of Definition \ref{def:Regular} below. The differential in space $\mu \in \Pcal_c(\R^d) \mapsto \D_x v[\mu](t,x)$ is narrowly continuous for $\Lcal^1$-almost every $t \in [0,T]$ and all $x \in \R^d$. \smallskip
\item[] \textnormal{\textbf{(B)}} : The barycenter $x \mapsto \bar{\Bgamma}^{\circ}_{\varphi}(x)$ of the minimal selection $\Bpartial^{\circ}\varphi(\mu)$ in the extended subdifferential of the terminal cost $\varphi(\cdot)$ at some  measure $\mu \in \Pcal_c(\R^d)$ is continuous. \\
The barycenter $x \mapsto \bar{\Bgamma}^{\circ}_L(x)$ of the minimal selection $\Bpartial_{\mu}^{\circ} L(\mu,\omega)$ in the extended subdifferential of the running cost $L(\cdot,\omega)$ at some $\mu \in \Pcal_c(\R^d)$ is continuous. \\
The barycenters $(x,y) \mapsto \bar{\Bgamma}_{(t,x)}^{i,\circ}(y)$ of the minimal selections $\Bpartial^{\circ}_{\mu}v^i[\mu](t,x)$ in the extended subdifferentials of the components $v^i$ define continuous mappings for $\Lcal^1$-almost every $t \in [0,T]$. \smallskip
\item[] \textnormal{\textbf{(D)}} : The maps $\mu \in \Pcal_c(\R^d) \mapsto \varphi(\mu)$, $\mu \in \Pcal_c(\R^d) \mapsto L(\mu,\omega)$ and $\mu \in \Pcal_c(\R^d) \mapsto v[\mu](t,x)$ are differentiable along measure curves generated by Lipschitz-in-time, continuous and bounded perturbations of the identity for $\Lcal^1 \times \mu$-almost every $(t,x) \in [0,T] \times \R^d$ and any $\omega \in U$, i.e.
\begin{equation*}
\begin{aligned}
\frac{d^+}{d \epsilon} \left[ \phi(\G(\epsilon,\cdot)_{\#} \mu) \right] & = \frac{d^-}{d \epsilon} \left[ \phi(\G(\epsilon,\cdot)_{\#} \mu) \right] \\ 
\frac{d^+}{d \epsilon} \left[ L(\G(\epsilon,\cdot)_{\#} \mu,\omega) \right] & = \frac{d^-}{d \epsilon} \left[ L(\G(\epsilon,\cdot)_{\#} \mu,\omega) \right] 
\end{aligned}
\end{equation*} 
and 
\begin{equation*}
\frac{d^+}{d \epsilon} \left[ v[\G(\epsilon,\cdot)_{\#} \mu](t,x) \right] = \frac{d^-}{d \epsilon} \left[ v[\G(\epsilon,\cdot)_{\#} \mu](t,x) \right],
\end{equation*}
whenever $(\G(\epsilon,\cdot))_{(-\bar{\epsilon},\bar{\epsilon})}$ is a Lipschitz family of continuous and bounded maps, differentiable at $\epsilon = 0$ and such that $\G(0,\cdot) = \Id$.
\end{enumerate}
\end{framed}

Then, there exist a constant $R > 0$ depending on $\mu^0$, $T$, $\U$, $v[\cdot](\cdot,\cdot)$, $\varphi(\cdot)$, $L(\cdot,\cdot)$ and  a curve $\nu^*(\cdot) \in \Lip([0,T],\Pcal(\overline{B_{2d}(0,R)})$ Lipschitzian with respect to the $W_1$-metric satisfying the following conditions : 
\begin{framed}
\begin{enumerate}
\item[(i)] It solves the forward-backward system of continuity equations
\begin{equation} \label{eq:Thm_HamiltonianODE1}
\left\{
\begin{aligned}
& \partial_t \nu^*(t) + \nabla_{(x,r)} \cdot \left( \J_{2d} \tilde{\nabla}_{\nu} \H_c(t,\nu^*(t),u^*(t))\nu^*(t) \right) = 0 ~ \text{in $[0,T] \times \R^{2d}$}, \\
& \pi^1_{\#} \nu^*(0) = \mu^0, \\
& \pi^2_{\#} \nu^*(T) = (-\bar{\Bgamma}^{\circ}_{\varphi})_{\#} \mu^*(T),
\end{aligned}
\right.
\end{equation}
where $\bar{\Bgamma}^{\circ}_{\varphi}(\cdot)$ is the barycenter of the minimal selection $\Bpartial^{\circ} \varphi(\mu^*(T))$ of the final cost $\varphi(\cdot)$ at $\mu^*(T)$ and $\J_{2d}$ is the symplectic matrix in $\R^{2d}$. 

\smallskip

The \textnormal{compactified Hamiltonian} of the system $\H_c(\cdot,\cdot,\cdot)$ is defined by 
\begin{equation} 
\label{eq:General_Hamiltonian_Compact}
\hspace{-0.45cm} \H_c (t,\nu,\omega) = ~  
\left\{
\begin{aligned} 
& \H(t,\nu,\omega) ~ \text{if $\nu \in \Pcal(\overline{B_{2d}(0,R)})$}, \\
& +\infty \hspace{2.4cm} \text{otherwise},
\end{aligned}
\right. 
\end{equation}
for any $(t,\nu,\omega) \in [0,T] \times \Pcal_c(\R^{2d}) \times U$ where 
\begin{equation}
\label{eq:General_Hamiltonian}
\H(t,\nu,\omega) = \INTDom{\langle r , v[\pi^1_{\#}\nu](t,x) + \omega(x) \rangle}{\R^{2d}}{\nu(x,r)} - L(\pi^1_{\#}\nu,\omega)
\end{equation}
is the \textnormal{infinite dimensional Hamiltonian} of the system for any $(t,\nu,\omega) \in [0,T] \times \Pcal_c(\R^{2d}) \times U$.

The vector field $\tilde{\nabla}_{\nu} \H_c(t,\nu^*(t),u^*(t))(\cdot,\cdot)$ is defined by 
\begin{equation*}
\begin{aligned}
& \tilde{\nabla}_{\nu} \H_c(t,\nu^*(t),u^*(t)) : (x,r) \in \supp(\nu^*(t)) \mapsto \\
& \begin{pmatrix} 
\D_x u^*(t,x)^{\top} r + \D_x v[\pi^1_{\#}\nu^*(t)](t,x)^{\top} r + \BGamma^{\circ}_v[\nu^*(t)](t,x) - \bar{\Bgamma}^{\circ}_L(t,x) \\ v[\pi^1_{\#} \nu^*(t)](t,x) + u^*(t,x)
\end{pmatrix}
\end{aligned}
\end{equation*}
where $t \in [0,T] \mapsto \bar{\Bgamma}_L^{\circ}(t,\cdot)$ is a measurable selection of the barycenters of $\Bpartial^{\circ}_{\mu} L(\mu^*(t),u^*(t))$. 

The map $\BGamma_v^{\circ}[\nu](\cdot,\cdot)$ is defined for any $\nu \in  \Pcal(\overline{B_{2d}(0,R)})$ by
\begin{equation*}
\BGamma^{\circ}_v [\nu] : (t,x) \in [0,T] \times \pi^1(\overline{B(0,R)}) \mapsto \INTDom{ \left( \BGamma^{\circ}_{(t,y)}(x) \right)^{\top} p \, }{\R^{2d}}{\nu(y,p)}
\end{equation*}
where for $\Lcal^1 \times \pi^1_{\#} \nu$-almost every $(t,y) \in [0,T] \times \pi^1(\overline{B_{2d}(0,R)})$ we define $\BGamma^{\circ}_{(t,y)} : x \in \supp(\pi^1_{\#}\nu(t)) \mapsto (\bar{\Bgamma}^{i,\circ}_{(t,y)}(x))_{1 \leq i \leq d}$ as the matrix-valued map made of the barycenters of the minimal selections $\Bpartial_{\mu}^{\circ} v^i[\pi^1_{\#}\nu](t,y)$ in the extended subdifferentials of the components $(v^i)$ of the non-local velocity field.

\medskip

\item[(ii)] It satisfies the Pontryagin maximization condition
\begin{equation} \label{eq:General_Maximization}
\H_c(t,\nu^*(t),u^*(t)) = \max\limits_{\omega \in U} \left[ \H_c(t,\nu^*(t),\omega) \right],
\end{equation}
for $\Lcal^1$-almost every $t \in [0,T]$.
\end{enumerate}
\end{framed}
\end{thm}

The general hypotheses \textbf{(H)} are rather cumbersome and can sometimes be hard to verify. Nevertheless, they are satisfied by a good score of functionals of great interest in various application fields. We present some relevant examples in Section \ref{section:Examples}.

\begin{rmk}[On the smoothness assumption \textnormal{\textbf{(U)}}]
The reason why we chose to impose the strong $C^{1,1}$-smoothness assumption on the set of admissible controls is twofold.

First, the main scope of this paper is to provide first-order optimality conditions for infinite-dimensional problems arising as mean-field limits of finite dimensional systems. Even though very general existence results \`a la DiPerna-Lions-Ambrosio \cite{AmbrosioPDE,DiPernaLions} are available for Cauchy problems of the form \eqref{eq:Nonlocal_Cauchy_problem}, they only deal with macroscopic quantities which are related to the underlying microscopic ones only for almost every curve in a suitable space of curves. The desired exact micro-macro correspondence which we aim at preserving can only hold in the presence of Cauchy-Lipschitz smoothness assumptions on the driving vector fields, see \cite{AmbrosioPDE}.

Second, the classical geometric proof of the maximum principle consisting in performing local-in-time perturbations of an optimal trajectories can only be carried out under $C^1$-regularity assumptions, due to the  non-linearity of the problem studied here. Even though the derivation of a maximum principle under a merely Lipschitz-regularity assumption on the optimal control in the spirit of the non-smooth maximum principle (see e.g. \cite{Clarke}) might be available in this context, it would require a completely different approach and much more technical arguments.

Let it be noted that these assumptions are verified in the classical setting of systems that are linear or affine with respect to the controls, i.e. where the controlled term is of the form $u : (t,x) \mapsto \sum_{k=1}^m u_k(t) F_k(x)$ where the $(F_k(\cdot))_{1 \leq k \leq m}$ are $C^1$ vector fields.
\end{rmk}

\begin{rmk}[Almost-Hamiltonian flow]
Observe that in our formulation of the PMP, the vector field $\tilde{\nabla}_{\nu} \H_c(t,\nu^*(t),u^*(t))$ is not the Wasserstein gradient of the compactified Hamiltonian $\H_c(\nu^*(t),u^*(t))$, since in general the barycenter of a minimal selection is not in the classical subdifferential. However, in any context where the minimal selections of the cost and dynamics functionals are induced by maps, which will automatically be their barycenters, or when they are strong subdifferentials, it can be shown by standard methods that $\tilde{\nabla}_{\nu} \H_c(t,\nu^*(t),u^*(t))$ is in fact the Wasserstein gradient of the compactified Hamiltonian at $(t,\nu^*(t),u^*(t))$ for $\Lcal^1$-almost every $t\in [0,T]$. 
\end{rmk}

We first describe in Section \ref{subsection:SimplePMP} our scheme of proof on a simplified problem $(\pazocal{P}_1)$ where there are no interaction field $v[\cdot](\cdot,\cdot)$ and no running cost $L(\cdot,\cdot)$. We then proceed to prove the PMP for the more general problem $(\pazocal{P})$ in Section \ref{subsection:GeneralPMP}. In what follows, we shall restrict our attention to the Wasserstein space $\Pcal_c(\R^d)$ endowed with the $W_1$-metric.

%%%%%%%%%%%%%%%%%%%%%%%%%%%%%%%%%%%%%%%%%%%%%%%%%%%%%%%%%%%%%%%%%%%%%%%%%%%%

\subsection{The Pontryagin Maximum Principle with no interaction field and no running cost}
\label{subsection:SimplePMP}

We start by proving the Pontryagin Maximum Principle for a simplified version of the optimal control problem $(\pazocal{P})$ presented in the introduction. We consider the following optimal control problem in the space of probability measures
\begin{equation} \label{eq:Simple_OCP}
(\pazocal{P}_1) ~~ \left\{
\begin{aligned}
& \underset{u \in \U}{\text{min}} \left[ \varphi(\mu(T)) \right], \\
& \text{s.t.} 
\left\{ 
\begin{aligned}
& \partial_t \mu(t) + \nabla \cdot (u(t,\cdot) \mu(t)) = 0, \\
& \mu(0) = \mu^0 \in \mathcal{P}_c(\R^d),
\end{aligned}
\right.
\end{aligned}
\right.
\end{equation}
and show that the Pontryagin-type optimality conditions provided in the following theorem hold.

\begin{thm}[Pontryagin Maximum Principle for $(\pazocal{P}_1)$]
\label{thm:Simple_PMP}
\hfill \\ Let $(u^*(\cdot),\mu^*(\cdot)) \in \U \times \Lip([0,T],\Pcal_c(\R^d))$ be an optimal pair control-trajectory for $(\pazocal{P}_1)$ and assume that hypotheses \textnormal{\textbf{(U),(C),(B)}} hold. Then, there exists a constant $R > 0$ and a curve $\nu^*(\cdot) \in \Lip([0,T],\Pcal(\overline{B_{2d}(0,R)}))$ satisfying the following statements :
\begin{enumerate}
\item[(i)] It solves the forward-backward system of continuity equations
\begin{equation} \label{eq:Thm_HamiltonianODE2}
\left\{
\begin{aligned}
& \partial_t \nu^*(t) + \nabla_{(x,r)} \cdot \left( \J_{2d} \nabla_{\nu} \H_c(\nu^*(t),u^*(t)) \nu^*(t) \right) = 0, ~~ & \text{in $[0,T] \times \R^{2d}$} \\
& \pi^1_{\#} \nu^*(0) = \mu^0, \\
& \pi^2_{\#} \nu^*(T) = (-\bar{\Bgamma}^{\circ}_{\varphi})_{\#} \mu^*(T),
\end{aligned}
\right.
\end{equation}
where $\J_{2d}$ is the symplectic matrix of $\R^{2d}$.

The \textnormal{compactified Hamiltonian} $\H_c(\cdot,\cdot)$ of the system is defined by 
\begin{equation}
\label{eq:Simple_Compact_Hamiltonian}
\H_c : (\nu,\omega) \in \Pcal_c(\R^{2d}) \times U \mapsto ~ 
\left\{
\begin{aligned} 
& \H(\nu,\omega) ~~ & \text{if $\nu \in \Pcal(\overline{B_{2d}(0,R)})$} \\
& +\infty ~~ & \text{otherwise}
\end{aligned}
\right. ~,
\end{equation}
where 
\begin{equation} \label{eq:Simple_Hamiltonian}
\H : (\nu,\omega) \in \Pcal_c(\R^{2d}) \times U \mapsto \INTDom{\langle r , \omega(x) \rangle}{\R^{2d}}{\nu(x,r)}
\end{equation}
is the \textnormal{infinite dimensional Hamiltonian} of the system.

The vector field $(x,r) \mapsto \nabla_{\nu} \H_c(\nu^*(t),u^*(t))(x,r) = (\D_x u^*(t,x)^{\top}r,u^*(t,x))$ is the \textnormal{Wasserstein gradient} of the compactified Hamiltonian for $\Lcal^1$-almost every $t \in [0,T]$, i.e. $\Bpartial^{\circ}_{\nu} \H_c(\nu^*(t),u^*(t)) = (I_{2d} \times \nabla_{\nu} \H_c(\nu^*(t),u^*(t)))_{\#} \nu^*(t)$.

\item[(ii)] It satisfies the Pontryagin maximization condition 
\begin{equation} \label{eq:Simple_Maximization}
\H_c(\nu^*(t),u^*(t)) = \max\limits_{\omega \in U} \left[ \H_c(\nu^*(t),\omega) \right]
\end{equation}
for $\Lcal^1$-almost every $t \in [0,T]$.
\end{enumerate}
\end{thm}

We split the proof of this result into several steps. In Step 1, we introduce the concept of \textit{needle-like variation} of an optimal control and compute explicitly the corresponding family of perturbed measures. In Step 2 we study the first order perturbation of the final cost induced by the needle-like variation. We introduce in Step 3 a suitable costate propagating this information backward to the base point of the needle-variation. In Step 4, we show that the curve introduced in Step 3 satisfies the conditions $(i)$ and $(ii)$ of the PMP. 

\begin{flushleft}
\underline{\textbf{Step 1 : Needle-like variations}}
\end{flushleft}

We start by considering an optimal pair control-trajectory $(u^*(\cdot),\mu^*(\cdot)) \in \U \times \Lip([0,T],\mathcal{P}_{\text{c}}(\R^d))$ along with the constant $R_T > 0$ such that $\supp(\mu(t)) \subset \overline{B(0,R_T)}$ for all times $t \in [0,T]$. Fix a control $\omega \in U$, a Lebesgue point $\tau \in [0,T]$ of $t \mapsto u^*(t) \in C^1(\R^d,\R^d)$ in the Bochner sense (see, e.g. \cite{DiestelUhl}) and a parameter $\epsilon \in [0,\bar{\epsilon})$ with $\bar{\epsilon} > 0$ small. We define the \textit{needle-like variation} of parameters $(\omega,\tau,\epsilon)$ of $u^*$ as follows
\begin{equation} \label{eq:Needle_variation}
\tilde{u}^{\omega,\tau}_{\epsilon} \equiv \tilde{u}_{\epsilon} : t \mapsto 
\left\{
\begin{aligned}
& \omega ~~ & \text{if $t \in [\tau-\epsilon,\tau]$}, \\
& u^*(t) ~~ & \text{otherwise}.
\end{aligned}
\right.
\end{equation}

We denote by $t \mapsto \tilde{\mu}_t(\epsilon)$ the corresponding solution of the continuity equation starting from $\mu^0$ at time $t = 0$. Notice that $\tilde{u}_{\epsilon}(\cdot) \in L^1([0,T],C^1(\R^d,\R^d) \cap L^{\infty}(\R^d,\R^d))$, thus the corresponding continuity equation is still well-posed.

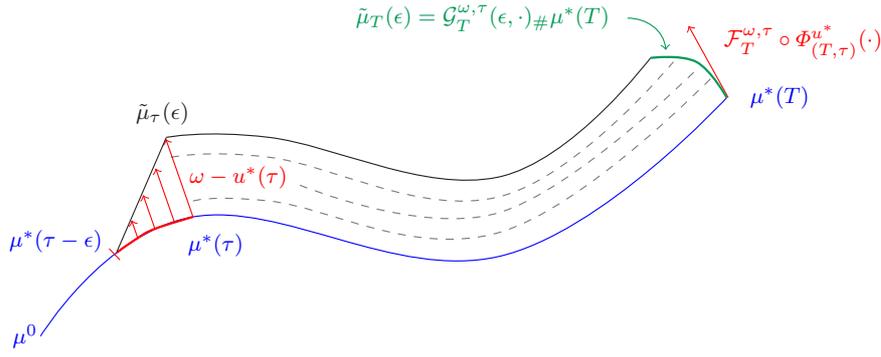
\begin{figure}[!h]
\centering
\resizebox{0.8 \textwidth}{.3\textwidth}{
\begin{tikzpicture}
\draw [blue] plot [smooth, tension=0.7] coordinates { (-1,-1)(1,0.5) (5,0)(8,2)};
\draw [red,thick] plot [smooth, tension=0.7] coordinates { (0,0.05)(0.25,0.22) (0.5,0.35) (1,0.5)};
\draw[red] (-0.1,0.1)--+(-45:0.2);
\draw[blue] (-1.2,-1) node{\footnotesize $\mu^0$};
\draw[blue] (-0.8,0.2) node{\footnotesize $\mu^*(\tau-\epsilon)$};
\draw[blue] (1.3,0.15) node{\footnotesize $\mu^*(\tau)$};
\draw[blue] (8.2,2) node[right] {\footnotesize $\mu^*(T)$};
\draw[red,->] (0.28,0.22)--+(110:0.25);
\draw[red,->] (0.5,0.35)--+(110:0.45);
\draw[red,->] (0.75,0.45)--+(110:0.7);
\draw[red,->] (1,0.5)--+(110:1.04);
\draw[red] (1.6,1) node {\footnotesize $\omega - u^*(\tau)$};
\draw[red,->] (8.02,2)--+(120:1.04);
\draw[red] (9,2.7) node {\footnotesize $\F^{\omega,\tau}_T \circ \Phi^{u^*}_{(T,\tau)}(\cdot)$};
\draw[Black] (0.6,1.8) node {\footnotesize $\tilde{\mu}_{\tau}(\epsilon)$};
\draw[Black](-0.,0.06)--(0.65,1.5);
\draw [Black] plot [smooth, tension=0.7] coordinates { (0.65,1.5) (2,1.5) (5,1)(7,2.5)};
\draw [Black , opacity = 0.6 , dashed] plot [smooth, tension=0.7] coordinates { (0.71,1.25) (2,1.25) (5,0.75)(7.3,2.45)};
\draw [Black , opacity = 0.6 , dashed] plot [smooth, tension=0.7] coordinates { (2.4,0.9) (5,0.55)(7.6,2.35)};
\draw [Black , opacity = 0.6 , dashed] plot [smooth, tension=0.7] coordinates { (1,0.7) (2,0.7) (5,0.3)(7.8,2.25)};
\draw[ForestGreen] (4.8,3) node {\footnotesize $\tilde{\mu}_T(\epsilon) = \G^{\omega,\tau}_T(\epsilon,\cdot)_{\#} \mu^*(T)$};
\draw[-> , ForestGreen](6.7,3) to [out = 0 , in = 95](7.2,2.6);
\draw [ForestGreen,thick] plot [smooth, tension=0.7] coordinates { (8,2) (7.6,2.45) (7,2.5)};
\end{tikzpicture}
}
\caption{\textit{Illustration of the effect of a needle-like variation on a measure curve.}}
\end{figure}

The link between the perturbed measure $\tilde{\mu}_T(\epsilon)$ and the optimal measure $\mu^*(T)$ at time $T$ is given in the following Lemma.

\begin{lem} \label{lem:Perturbed_final_state}
There exists a family of functions $\G^{\omega,\tau}_T(\cdot,\cdot) \in \Lip((-\bar{\epsilon},\bar{\epsilon}),C^0(\R^d,\R^d))$ such that they  are $C^1$-diffeomorphisms over $\overline{B(0,R_T)}$ for all $\epsilon \geq 0$ and it holds
\begin{equation}
\label{eq:Perturbed_final_state}
\tilde{\mu}_T(\epsilon) = \G^{\omega,\tau}_T(\epsilon,\cdot)_{\#}\mu^*(T).
\end{equation}
Moreover, there exists a constant $R_T^{\Phi} > 0$ depending on $R_T,L_U$ and $v[\cdot](\cdot,\cdot)$ such that for all $\epsilon \in (-\bar{\epsilon},\bar{\epsilon})$ it holds $\supp(\G^{\omega,\tau}_T(\epsilon,\cdot)_{\#}\mu^*(T)) \subset \overline{B(0,R_T^{\Phi})}$. 

This family of maps satisfies the following Taylor expansion with respect to the $L^2(\R^d,\R^d;\mu^*(T))$-norm
\begin{equation*}
\G^{\omega,\tau}_{T}(\epsilon,\cdot) = \Id + \epsilon \F^{\omega,\tau}_T \circ \Phi^{u^*}_{(T,\tau)}(\cdot) + o(\epsilon),
\end{equation*}
where
\begin{equation}
\label{eq:FirstOrder_Needle_Expression}
\F^{\omega,\tau}_T : x \in \supp(\mu^*(\tau)) \mapsto \D_x \Phi^{u^*}_{(\tau,T)}(x) \cdot \left[ \omega(x) - u^*(\tau,x) \right] 
\end{equation}
is a $C^0$ mapping.
\end{lem}

\begin{proof}
 
By definition of $\tilde{u}_{\epsilon}(\cdot,\cdot)$ in \eqref{eq:Needle_variation}, it holds that 
\begin{equation*}
\tilde{\mu}_T(\epsilon) = \Phi^{u^*}_{(\tau,T)} \circ \Phi^{\omega}_{(\tau-\epsilon,\tau)} \circ \Phi^{u^*}_{(\tau,\tau-\epsilon)} \circ \Phi^{u^*}_{(T,\tau)}(\cdot)_{\#} \mu^*(T) ~ \text{for all $\epsilon \in [0,\bar{\epsilon})$}.
\end{equation*}
Thus, by choosing $ \G^{\omega,\tau}_T(\epsilon,\cdot) = \Phi^{u^*}_{(\tau,T)} \circ \Phi^{\omega}_{(\tau-\epsilon,\tau)} \circ \Phi^{u^*}_{\tau-\epsilon} \circ \Phi^{u^*}_{(T,\tau)}(\cdot)$, formula \eqref{eq:Perturbed_final_state} holds true for $\epsilon \in [0,\bar{\epsilon})$. Moreover, since the definition of $\epsilon \mapsto \G^{\omega,\tau}_T(\epsilon,\cdot)$ only involves functions that are continuous and uniformly bounded over $\overline{B(0,R_T)}$, the perturbed measures $\tilde{\mu}_T(\cdot)$ are compactly supported in some bigger ball $\overline{B(0,R'_T)}$ for all $\epsilon \in [0,\bar{\epsilon})$ as well.

Recalling the definition of the flow $x \mapsto \Phi^v_{(0,t)}(x)$, one has by Lebesgue's Differentiation Theorem (see e.g. \cite[Chapter 1.7]{EvansGariepy}) :
\begin{equation*}
\left\{
\begin{aligned}
&\Phi^{u^*}_{(\tau,\tau-\epsilon)}(x) = x - \INTSeg{u^* \left( t,\Phi^{u^*}_{(\tau,t)}(x) \right)}{t}{\tau-\epsilon}{\tau} = x - \epsilon u^*(\tau,x) + o(\epsilon) \\
& \Phi^{\omega}_{(\tau-\epsilon,\tau)}(x) = x + \INTSeg{\omega \left( \Phi_{(\tau-\epsilon,t)}^{\omega}(x) \right)}{t}{\tau-\epsilon}{\tau} = x + \epsilon \omega(x) + o(\epsilon)  
\end{aligned}
\right.
\end{equation*}
since $t \mapsto \Phi_{(\tau-\epsilon,t)}^{u^*}(x)$ and $t \mapsto \Phi_{(t,\tau)}^{u^*}(x)$  are $C^0$ for any $x \in \overline{B(0,R_T)}$ and $\tau \in [0,T]$ Lebesgue point of $t \mapsto u^*(t) \in C^1(\R^d,\R^d)$ in the Bochner sense. Chaining these two expansions and recalling that $\omega(\cdot)$ and $\Phi^{u^*}_{(\tau,T)}(\cdot)$ are $C^1$-smooth yields
\begin{equation*}
\Phi^{u^*}_{(\tau,T)} \circ \Phi^{\omega}_{(\tau-\epsilon,\tau)} \circ \Phi^{u^*}_{(\tau,\tau-\epsilon)}(x) = \Phi^{u^*}_{(\tau,T)}(x) + \epsilon \D_x \Phi^{u^*}_{(\tau,T)}(x) \cdot \left[ \omega(x) - u^*(\tau,x) \right] + o(\epsilon).
\end{equation*}
Thus,
\begin{equation*}
\G^{\omega,\tau}_T(\epsilon,x) = x + \epsilon \F^{\omega,\tau}_T \circ  \Phi^{u^*}_{(T,\tau)}(x) + o(\epsilon) ~~ \text{for any $x \in \supp(\mu^*(T))$},
\end{equation*}
where we choose 
\begin{equation*}
\F^{\omega,\tau}_T : x \mapsto \D_x \Phi^{u^*}_{(\tau,T)}(x) \cdot \left[ \omega(x) - u^*(\tau,x)\right].
\end{equation*}

We can now extend $\G_T^{\omega,\tau}(\cdot,\cdot)$ from $[0,\bar{\epsilon})$ to $(-\bar{\epsilon},\bar{\epsilon})$ in such a way that the left and right derivatives at $\epsilon = 0$ coincide, by defining e.g.
\begin{equation*}
\G_T^{\omega,\tau}(\epsilon,\cdot) \equiv \Id + \epsilon \F^{\omega,\tau}_T \circ \Phi^{u^*}_{(T,\tau)}(\cdot)
\end{equation*}
whenever $\epsilon \in (-\bar{\epsilon},0]$. Notice that since $\F^{\omega,\tau}_T \circ \Phi^{u^*}_{(T,\tau)}(\cdot)$ is $C^0$, both $x \in \overline{B(0,R_T)}\mapsto \G^{\omega,\tau}_T(\epsilon,x)$ and $x \in \overline{B(0,R_T)}\mapsto \derv{}{\epsilon}{}[\G^{\omega,\tau}_T(\epsilon,x)]_{\epsilon = 0}$ define $C^0$ mappings for all $\epsilon \in (-\bar{\epsilon},\bar{\epsilon})$. Moreover, the continuity and uniform boundedness of $\D_x \Phi^{u^*}_{(\tau,T)}(\cdot)$ over $\overline{B(0,R_T)}$ along with hypothesis \textbf{(U)} imply that $\F^{\omega,\tau}_T(\cdot)$ is bounded. Hence, there exists a constant $R_T^{\Phi} > 0$ such that $\supp(\G^{\omega,\tau}_T(\epsilon,\cdot)_{\#} \mu^*(T))$ $\subset \overline{B(0,R_T^{\Phi})}$ for all $\epsilon \in (-\bar{\epsilon},\bar{\epsilon})$. Moreover, the fact that $\G(\epsilon,\cdot)$ and $\F^{\omega,\tau}_T(\cdot)$ are continuous and bounded yields that they are uniformly integrable with respect to the compactly supported measure $\mu^*(T)$. An application of Proposition \ref{prop:Convergence_results}$-(iii)$ allows to conclude that this expansion holds in $L^2(\R^d,\R^d;\mu^*(T))$, which achieves the proof. 
\end{proof}

We end this first step by a Lemma which is a direct consequence of Proposition \ref{prop:Differential_local_flow}. 

\begin{lem} \label{lem:Differential_F}
For any $x \in \supp(\mu^*(\tau))$, the trajectory $t \mapsto \F^{\omega,\tau}_t(x)$ is the unique solution of the Cauchy problem
\begin{equation}
\label{eq:Differential_F}
\partial_t \F^{\omega,\tau}_t(x) = \D_x u^* \left( t,\Phi_{(\tau,t)}^{u^*}(x) \right) \F^{\omega,\tau}_t(x) ~,~ \F^{\omega,\tau}_{\tau}(x) = \omega(x) - u^*(\tau,x).
\end{equation}
\end{lem}

\begin{proof}
It is sufficient to apply Proposition \ref{prop:Differential_local_flow} and to remark that here $v \equiv u^*$.
\end{proof}

\begin{flushleft}
\underline{\textbf{Step 2 : First-order optimality condition}}
\end{flushleft}

Thanks to the optimality of $u^*(\cdot)$, for each $\epsilon \in (0,\bar{\epsilon})$ it holds
\begin{equation} \label{eq:Optimality_endcost}
\frac{\varphi(\tilde{\mu}_T(\epsilon)) - \varphi(\mu^*(T))}{\epsilon} \geq 0,
\end{equation}
where $\epsilon \in (0,\bar{\epsilon}) \mapsto \tilde{\mu}_T(\epsilon) = \G^{\omega,\tau}_T(\epsilon,\cdot)_{\#} \mu^*(T)$.

Recalling that the measures $\epsilon \mapsto \tilde{\mu}_T(\epsilon)$ are uniformly compactly supported, that $\varphi(\cdot)$ satisfies hypotheses \textbf{(C)} and that the map $\epsilon \in (-\bar{\epsilon},\bar{\epsilon}) \mapsto \varphi(\tilde{\mu}_T(\epsilon))$ is differentiable at $\epsilon =0$ by hypothesis \textbf{(D)}, we can apply the chainrule given in Proposition \ref{prop:Chainrule} to the endpoint cost : 
\begin{equation} \label{eq:Gradient_terminal_cost}
0 ~ \leq ~ \lim\limits_{\epsilon \downarrow 0} \left[ \frac{\varphi(\tilde{\mu}_T(\epsilon)) - \varphi(\mu^*(T))}{\epsilon} \right] = \INTDom{\langle \bar{\Bgamma}^{\circ}_{\varphi}(x) , \F^{\omega,\tau}_T \circ \Phi^{u^*}_{(T,\tau)}(x) \rangle}{\R^d}{\mu^*(T)(x)},
\end{equation}
where $\bar{\Bgamma}^{\circ}_{\varphi} \in L^2(\R^d,\R^d;\mu^*(T))$ is the barycenter of the minimal selection $\Bpartial^{\circ} \varphi(\mu^*(T))$ in the extended subdifferential of $\varphi(\cdot)$ at $\mu^*(T)$. 

We recover a formula similar to the classical finite dimensional case. The next step is to introduce a suitable \textit{costate} along with its backward dynamics that will \textit{propagate} this first-order information to the base-point $\tau$ of the needle-like variation while generating a Hamiltonian-like dynamical structure. 

\begin{flushleft}
\underline{\textbf{Step 3 : Backward dynamics and Pontryagin maximization condition}}
\end{flushleft}

Equation \eqref{eq:Gradient_terminal_cost} provides us with a first-order optimality condition which involves all the needle parameters $(\omega,\tau) \in U \times [0,T]$. We will show that it implies, along with the choice of a suitable costate, the maximization condition \eqref{eq:Simple_Maximization}. 

To this aim, we build a curve $\nu^* \in \Lip([0,T],\Pcal_c(\R^{2d}))$ solution of the forward-backward system of continuity equations
\begin{equation} 
\label{eq:SimplePMP_HamiltonianDynamics}
\left\{
\begin{aligned}
& \partial_t \nu^*(t) + \nabla \cdot (\V^*(t,\cdot,\cdot)) \nu^*(t)) = 0 ~~ & \text{in $[0,T] \times \R^{2d}$}, \\
& \pi^1_{\#} \nu^*(t) = \mu^*(t) ~~ & \text{for all $t \in [0,T]$}, \\
& \nu^*(T) = (\Id \times (-\bar{\Bgamma}^{\circ}_{\varphi}))_{\#} \mu^*(T),
\end{aligned}
\right.
\end{equation}
associated to the vector field 
\begin{equation}
\V^* : (t,x,r) \in [0,T] \times \R^{2d} \mapsto (u^*(t,x),-\D_x u^*(t,x)^{\top}r).
\label{e-Vstar}
\end{equation} Notice that, contrarily to system \eqref{eq:Thm_HamiltonianODE2}, we impose the more restrictive product structure on the terminal datum.

This system is peculiar in the sense that the driving vector field $\V^*(\cdot,\cdot,\cdot)$ does not satisfy verbatim the hypotheses \textbf{(H')} of Theorem \ref{thm:Existence_uniqueness_non-local_PDE}. However, it exhibits a \textit{cascade} structure, in the sense that one can first determine uniquely $\mu^*(\cdot)$ and then build $\nu^*(\cdot)$ by disintegration. This fact is underlined by the condition $\pi^1_{\#} \nu^*(t) = \mu^*(t)$ for all times $t \in [0,T]$. We make this statement precise in the next Lemma. 

\begin{lem}[Definition and well-posedness of solutions of \eqref{eq:SimplePMP_HamiltonianDynamics}]
\label{lem:Wellposedness_Simple}
\hfill \\ Let $(u^*(\cdot),\mu^*(\cdot))$ be an optimal pair control-trajectory for $(\pazocal{P}_1)$. For $\mu^*(T)$-almost every $x \in \R^d$, we consider the family of backward flows $(\Psi^x_{(T,t)}(\cdot))_{t \leq T}$ associated to the Cauchy problems
\begin{equation}
\label{eq:Marginal_TransportEq}
\partial_t w_x(t,r)= -\D_x u^*(t,\Phi^{u^*}_{(T,t)}(x)) w_x(t,r) ~,~ w_x(T,r) = r,
\end{equation}
and define the associated curves of measures $\sigma^*_x : t \mapsto \Psi^x_{(T,t)}(\cdot)_{\#} \delta_{(-\bar{\Bgamma}^{\circ}_{\varphi}(x))}$.

Then, $\nu^* : t \mapsto (\Phi^{u^*}_{(T,t)} , \Id)_{\#} \nu^*_T(t)$ is the unique solution of \eqref{eq:SimplePMP_HamiltonianDynamics} with\footnote{Namely, $\nu_T^*(t)$ is defined as the $\mu^*(T)$-almost uniquely determined measure which has $\mu^*(T)$ as its first marginal and which disintegration is given by $\{ \sigma^*_x(t)\}_x$ (see Definition \ref{def:Barycenter}).} $\nu^*_T(t) = \int \sigma^*_x(t) \textnormal{d} \mu^*(T)(x) \in \Pcal_c(\R^{2d})$ for all times $t \in [0,T]$. Moreover, there exists two constants $R_T',L_T' > 0$ such that
\begin{equation*}
\supp(\nu^*(t)) \subset \overline{B_{2d}(0,R'_T)} ~~ \text{and} ~~ W_1(\nu^*(t),\nu^*(s)) \leq L_T' |t-s|
\end{equation*}  
for all $s,t \in [0,T]$.
\end{lem}

\begin{proof}
We recall that by hypothesis \textbf{(U)}, the elements of $\U$ are uniformly sublinear and Lipschitz in space for $\Lcal^1$-almost every times $t \in [0,T]$. We recall that by Theorem \ref{thm:Existence_uniqueness_non-local_PDE}, this implies the existence of a constant $R_T > 0$ depending on $\supp(\mu^0),T$ and $L_U$ such that $\supp(\mu^*(\cdot)) \subset \overline{B(0,R_T)}$. 

For $\mu^*(T)$-almost every $x \in \R^d$, the Cauchy problem \eqref{eq:Marginal_TransportEq} has a unique solution and the corresponding curves $t \mapsto \sigma^*_x(t)$ are uniquely determined. Moreover, the uniform Lipschitzianity of the elements of $\U$ implies that these curves are uniformly compactly supported and Lipschitz in the $W_1$-metric uniformly with respect to $x \in \supp(\mu^*(T))$ with constants $\tilde{R}_T,\tilde{L}_T$ depending on $L_U,T$ and $\varphi(\cdot)$.

We now define the curve $\nu^*(\cdot)$ as in the statement of Lemma \ref{lem:Wellposedness_Simple} above and show that it is a uniformly compactly supported and Lipschitz solution of the forward-backward system \eqref{eq:SimplePMP_HamiltonianDynamics}. The fact that there exists $R_T > 0$ depending on $R_T$ and $\tilde{R}_T$ such that $\nu^*(\cdot)$ is uniformly compactly supported in $\overline{B_{2d}(0,R_T)}$ is a direct consequence of its definition. The Lipschitzianity in the $W_1$-metric comes from the following computations. For any $\xi \in \Lip(\R^{2d},\R)$ with $\Lip(\xi,\R^{2d}) \leq 1$, it holds
\begin{equation*}
\begin{aligned}
\INTDom{\xi(x,r)}{\R^{2d}}{(\nu^*(t)-\nu^*(s))(x,r)} & =  \INTDom{ \INTDom{\xi(\Phi^{u^*}_{(T,t)}(x),r)}{\R^d}{\sigma^*_x(t)(r)}}{\R^d}{\mu^*(T)(x)} \\
& \hspace{0.45cm} - \INTDom{ \INTDom{\xi(\Phi^{u^*}_{(T,s)}(x),r)}{\R^d}{\sigma^*_x(s)(r)}}{\R^d}{\mu^*(T)(x)}\\
& \leq \INTDom{ \INTDom{| \Phi^{u^*}_{(T,t)}(x) - \Phi^{u^*}_{(T,s)}(x) |}{\R^d}{\sigma_x^*(t)(r)}}{\R^d}{\mu^*(T)(x)} \\
& \hspace{0.45cm} + \INTDom{\Lip(\Phi^{u^*}_{(T,s)},\overline{B(0,R_T)}) W_1(\sigma^*_x(s),\sigma^*_x(t))}{\R^d}{\mu^*(T)(x)} \\
& \leq  L'_T |t-s|
\end{aligned}
\end{equation*}
where $L_T' > 0$ is a uniform constant depending on the time and space Lipschitz constants of the flows of diffeomorphims $(\Phi^{u^*}_{(T,t)}(\cdot))_{t \in [0,T]}$ and on $L_U$. Taking the supremum over all the 1-Lipschitz functions $\xi(\cdot,\cdot)$ and using the Kantorovich-Rubinstein duality \eqref{eq:Kantorovich_duality} yields the Lipschitzianity of $\nu^*(\cdot)$ in the $W_1$-metric.

Finally, remark that for any $\xi \in C^{\infty}_c(\R^{2d},\R)$ it holds
\begin{equation*}
\begin{aligned}
\derv{}{t} \left[ \INTDom{\xi(x,r)}{\R^{2d}}{\nu^*(t)(x,r)} \right] & =  \derv{}{t} \left[ \INTDom{ \INTDom{\xi(\Phi^{u^*}_{(T,t)}(x),r)}{\R^d}{\sigma^*_x(t)(r)}}{\R^d}{\mu^*(T)(x)} \right] \\
& =  \INTDom{ \INTDom{ \langle \nabla_x \xi(\Phi^{u^*}_{(T,t)}(x),r) , u^*(t,\Phi^{u^*}_{(T,t)}(x) \rangle }{\R^d}{\sigma^*_x(t)(r)}}{\R^d}{\mu^*(T)(x)} \\
& \hspace{0.45cm}  + \INTDom{ \INTDom{ \langle \nabla_r \xi(\Phi^{u^*}_{(T,t)}(x),r) , -\D_x u^*(t,\Phi^{u^*}_{(T,t)}(x))^{\top} r \rangle }{\R^d}{\sigma^*_x(t)(r)}}{\R^d}{\mu^*(T)(x)} \\
& = \INTDom{\left\langle \nabla_{(x,r)} \xi(x,r) , \begin{pmatrix} u^*(t,x)\\ -\D_x u^*(t,x)^{\top}r \end{pmatrix} \right\rangle}{\R^{2d}}{\nu^*(t)(x,r)}
\end{aligned}
\end{equation*}
which along with the fact that $\nu^*(T) = \nu^*_T(T) = (\Id \times (-\bar{\Bgamma}^{\circ}_{\varphi}))_{\#}\mu^*(T)$ achieves the proof.
\end{proof}

We now show that the curve of measures $\nu^*(\cdot)$ defined in Lemma \ref{lem:Wellposedness_Simple} is such that the map $\K_{\omega,\tau}(\cdot)$ defined by 
\begin{equation}
\label{eq:Simple_K}
\K_{\omega,\tau} : t \in [\tau,T] \mapsto \INTDom{\langle r , \F^{\omega,\tau}_t \circ \Phi^{u^*}_{(t,\tau)}(x) \rangle}{\R^{2d}}{\nu^*(t)(x,r)},
\end{equation}
is constant over $[\tau,T]$. We shall see in Step 4 that this is equivalent to the Pontryagin maximization condition \eqref{eq:Simple_Maximization}.

\begin{lem} \label{lem:Backward_VectorField_Simple}
The map $t \mapsto \K_{\omega,\tau}(t)$ defined in \eqref{eq:Simple_K} is constant over $[\tau,T]$ for any couple of needle parameters $(\omega,\tau)$.
\end{lem}
\begin{proof} 
Notice that by definition of $\nu^*(\cdot)$, the map $\K_{\omega,\tau}(\cdot)$ rewrites
\begin{equation}
\K_{\omega,\tau}(t) = \INTDom{\INTDom{\langle r , \F^{\omega,\tau}_t \circ \Phi^{u^*}_{(T,\tau)}(x) \rangle}{\R^d}{\sigma^*_x(t)(r)} }{\R^d}{\mu^*(T)(x)} ~~ \text{for all $t \in [\tau,T]$}.
\end{equation}
The maps $t \in [\tau,T] \mapsto \F^{\omega,\tau}_t \circ \Phi^{u^*}_{(T,\tau)}(x)$ and $t \in [\tau,T] \mapsto \sigma^*_x(t) $ are Lipschitz, uniformly with respect to $x \in \supp(\mu^*(T))$. The integrand $(x,r) \mapsto \langle r , \F^{\omega,\tau}_t \circ \Phi^{u^*}_{(T,\tau)}(x) \rangle$ is bounded with respect to $x$ and Lipschitz with respect to $r$, uniformly with respect to $t \in [\tau,T]$. Hence, $t \mapsto \K_{\omega,\tau}(t)$ is Lipschitz as well. It will therefore be constant provided that its derivative - which exists $\Lcal^1$-almost everywhere - is equal to zero. 

Observe that, using formula \eqref{eq:NonlocalPDE_distributions2} and the definition of $\V^*(\cdot,\cdot,\cdot)$ in \eqref{eq:Def_psi_epsilon}, it holds
\begin{equation} \label{eq:Time_derivative_K}
\begin{aligned}
\derv{}{t}\K_{\omega,\tau}(t) &  = \INTDom{ \INTDom{ \langle r , \partial_t \F^{\omega,\tau}_t \circ \Phi^{u^*}_{(T,\tau)}(x) \rangle}{\R^d}{\sigma^*_x(t)(r)} }{\R^d}{\mu^*(T)(x)} \\
& + \INTDom{ \INTDom{ \langle -\D_x u^*(t,\Phi^{u^*}_{(T,t)}(x))^{\top} r ,  \F^{\omega,\tau}_t \circ \Phi^{u^*}_{(T,\tau)}(x) \rangle}{\R^d}{\sigma^*_x(t)(r)} }{\R^d}{\mu^*(T)(x)}.
\end{aligned}
\end{equation} 

We recall the characterization of $\partial_t \F^{\omega,\tau}_t(\cdot)$ given in \eqref{eq:Differential_F} and plug it into \eqref{eq:Time_derivative_K}. This implies that $\derv{}{t}\K_{\omega,\tau}(t) = 0$ for $\Lcal^1$-almost every times $t \in [\tau,T]$, and thus that $\K_{\omega,\tau}(\cdot)$ is constant over $[\tau,T]$. 
\end{proof} 

\begin{flushleft}
\underline{\textbf{Step 4 : Proof of the Pontryagin Maximum Principle for $(\pazocal{P}_1)$}}
\end{flushleft}

We proved in Lemma \ref{lem:Wellposedness_Simple} the existence of a constant $R \equiv R_T' > 0$ such that the solution $\nu^*(\cdot)$ to \eqref{eq:SimplePMP_HamiltonianDynamics} satisfies $\supp(\nu^*(\cdot)) \subset \overline{B_{2d}(0,R)}$. We accordingly define the infinite dimensional Hamiltonian $\H : (\nu,\omega) \in \Pcal_c(\R^{2d}) \times U \mapsto \INTDom{\langle r , \omega(x) \rangle}{\R^{2d}}{\nu(x,r)}$ of the system and the compactified Hamiltonian $\H_c(\cdot,\cdot)$ by \eqref{eq:Simple_Compact_Hamiltonian}.

In Lemma \ref{lem:Backward_VectorField_Simple} we showed that, with this choice of forward-backward system \eqref{eq:SimplePMP_HamiltonianDynamics}, the map $\K_{\omega,\tau}(\cdot)$ defined in \eqref{eq:Simple_K} is constant over $[\tau,T]$ for any choice of $\omega \in U$ and $\tau \in [0,T]$ Lebesgue point of $u^*(\cdot)$. This implies in particular that $\K_{\omega,\tau}(\tau) = \K_{\omega,\tau}(T)$. Since we proved in \eqref{eq:Gradient_terminal_cost} that it holds
\begin{equation*}
0 \leq \INTDom{\langle \bar{\Bgamma}^{\circ}_{\varphi}(x) , \F^{\omega,\tau}_T \circ \Phi^{u^*}_{(T,\tau)}(x) \rangle}{\R^d}{\mu^*(T)(x)} = - K_{\omega,\tau}(T),
\end{equation*}
it directly follows that
\begin{equation*}
\K_{\omega,\tau}(\tau) = \INTDom{\langle r , \omega(x) - u^*(\tau,x) \rangle}{\R^{2d}}{\nu^*(\tau)(x,r)} \leq 0,
\end{equation*}
for all $\omega \in U$ and $\tau \in [0,T]$ Lebesgue point of $u^*(\cdot)$.

Recalling that $\Lcal^1$-almost $\tau \in [0,T]$ is a Lebesgue point in the Bochner sense for an $L^1$-function defined over the real line (see e.g. \cite[Chapter 2, Theorem 9]{DiestelUhl}), we recover the infinite dimensional maximization condition \eqref{eq:Simple_Maximization}
\begin{equation*}
\H_c(\nu^*(t),u^*(t)) = \max\limits_{\omega \in U} \left[ \H_c(\nu^*(t),\omega) \right] 
\end{equation*}
for $\Lcal^1$-almost every $t \in [0,T]$.

Invoking the $C^1$ regularity of the elements of $U$, it can be seen using Proposition \ref{prop:Subdifferential_Integral} that the minimal selection $\Bpartial^{\circ}_{\nu} \H_c(\nu^*(t),u^*(t))$ in the extended subdifferential of $\H_c(\cdot,u^*(t))$ exists at $\nu^*(t) \in \Pcal(\overline{B_{2d}(0,R)})$ for $\Lcal^1$-almost every $t \in [0,T]$ and that it is induced by the map
\begin{equation*}
\nabla_{\nu} \H_c(\nu^*(t),u^*(t)) : (x,r) \in \supp(\nu^*(t)) \mapsto \begin{pmatrix}
\D_x u^*(t,x)^{\top} r \\ u^*(t,x)
\end{pmatrix}.
\end{equation*}
Therefore, we recognize the Wasserstein Hamiltonian structure $\V^*(t,\cdot,\cdot) = \J_{2d} \nabla_{\nu} \H_c(\nu^*(t),u^*(t))(\cdot,\cdot)$ for $\Lcal^1$-almost every $t \in [0,T]$ where $\J_{2d}$ is the symplectic matrix in $\R^{2d}$. This ends our proof of Theorem \ref{thm:Simple_PMP}.

%%%%%%%%%%%%%%%%%%%%%%%%%%%%%%%%%%%%%%%%%%%%%%%%%%%%%%%%%%%%%%%%%%%%%%%%%%%%

\subsection{The general Pontryagin Maximum Principle}
\label{subsection:GeneralPMP}

After having exhibited the main mechanisms of our proof for the Pontryagin Maximum Principle for the simplified problem $(\pazocal{P}_1)$, we are ready to tackle the general case proposed in $(\pazocal{P})$. The latter is a generalization of $(\pazocal{P}_1)$ in the sense that we add a general running cost $L(\cdot,\cdot)$ and a general non-local interaction vector field $v[\cdot](\cdot,\cdot)$.

\begin{flushleft}
\textbf{\underline{Step 1 : Needle-like variations in the non-local case}}
\end{flushleft}

As in Section \ref{subsection:SimplePMP}, let us consider an optimal pair control-trajectory $(u^*(\cdot),\mu^*(\cdot))$, a Lebesgue point $\tau \in [0,T]$ of $u^*(\cdot)$ and an element $\omega \in U$. We introduce again the needle-like variation $\tilde{u}^{\omega,\tau}_{\epsilon}(\cdot)$ of $u^*(\cdot)$ with parameters $(\omega,\tau,\epsilon)$ for $\epsilon \in [0,\bar{\epsilon})$, as defined in \eqref{eq:Needle_variation}. Notice that this time, $\tau$ is a Lebesgue point for $t \mapsto v[\mu^*(t)](t,\cdot) + u^*(t,\cdot)$.

In keeping with the notations introduced in \eqref{eq:Flow_def} for flows associated to transport PDEs with non-local velocities, the family of perturbed measures $\epsilon \in [0,\bar{\epsilon}) \mapsto \tilde{\mu}_t(\epsilon)$ are defined for all times $t \in [\tau,T]$ by
\begin{equation*}
\tilde{\mu}_t(\epsilon) = \Phi^{v,u^*}_{(\tau,t)}[\tilde{\mu}_{\tau}(\epsilon)] \circ \Phi^{v,\omega}_{(\tau-\epsilon,\tau)}[\mu^*(\tau-\epsilon)] \circ \Phi^{v,u^*}_{(t,\tau-\epsilon)}[\mu^*(t)] (\cdot)_{\#} \mu^*(t).
\end{equation*}
One can readily check that under the sublinearity and regularity hypotheses imposed in \textbf{(U)} and \textbf{(F)}, there exists again a constant $\tilde{R}_T > 0$ such that $\supp(\tilde{\mu}_t(\epsilon)) \subset \overline{B(0,\tilde{R}_T)}$ for all $(t,\epsilon) \in [0,T] \times [0,\bar{\epsilon})$.

We now derive in Lemma \ref{lem:Perturbed_state_non-local} the perturbation stemming from the needle-like variation. We prove therein a result akin to Lemma \ref{lem:Differential_F} giving a precise ODE-type characterization of this perturbation. To do so, we use the results of Proposition \ref{prop:Directional_derivative_flow} concerning  the directional derivatives of the non-local flow combined to the classical result stated in Lemma \ref{lem:Differential_F} and the definition of needle-like variation. 

\begin{lem}[Perturbation induced by a needle-like variation in the non-local case]
\label{lem:Perturbed_state_non-local}
Let $(u^*(\cdot),\mu^*(\cdot))$ be an optimal pair control-trajectory for problem $(\pazocal{P})$ and $\tilde{u}_{\epsilon}(\cdot)$ be the needle-like perturbation of $u^*(\cdot)$ as introduced in \eqref{eq:Needle_variation}. 

Then, there exists for all times $t \in [\tau,T]$ a family of functions $\G^{\omega,\tau}_t(\cdot,\cdot) \in \Lip((-\bar{\epsilon},\bar{\epsilon}),C^0(\R^d,\R^d))$ such that they are $C^1$-diffeomorphisms over $\overline{B(0,R_T)}$ for all $\epsilon \geq 0$ and it holds
\begin{equation*}
\tilde{\mu}_t(\epsilon) = \G^{\omega,\tau}_t(\epsilon,\cdot)_{\#} \mu^*(t).
\end{equation*}
Besides, there exists a constant $R^{\Phi}_T > 0$ depending on $R_T,L_U$ and $v[\cdot](\cdot,\cdot)$ such that for all $(t,\epsilon) \in [\tau,T] \times (-\bar{\epsilon},\bar{\epsilon})$ it holds $\supp(\G^{\omega,\tau}_t(\epsilon,\cdot)_{\#} \mu^*(t)) \subset \overline{B(0,R_T^{\Phi})}$.

This family of maps satisfies the following Taylor expansion for all $t \in [\tau,T]$ with respect to the $L^2(\R^d,\R^d;\mu^*(t))$ norm : 
\begin{equation*}
\G^{\omega,\tau}_t(\epsilon,\cdot) = \Id + \epsilon \F^{\omega,\tau}_t \circ \Phi^{v,u^*}_{(t,\tau)}[\mu^*(t)](\cdot) + o(\epsilon), 
\end{equation*}
with
\begin{equation*}
\begin{aligned}
\F^{\omega,\tau}_t : x \in \supp(\mu^*(\tau)) & \mapsto \D_x \Phi^{v,u^*}_{(\tau,t)}[\mu^*(\tau)](x) \cdot \left[ \omega(x) - u^*(\tau,x) \right] + w_{\Phi}^{\omega,\tau}(t,x)
\end{aligned}
\end{equation*}
where $w_{\Phi}^{\omega,\tau}(t,x)$ is the derivative at $\epsilon = 0$ of the map $\epsilon \in (-\bar{\epsilon},\bar{\epsilon}) \mapsto \Phi^{v,u^*}_{(\tau,t)}[\tilde{\mu}_{\tau}(\epsilon)](x)$ as described in Proposition \ref{prop:Directional_derivative_flow}.

Moreover, the map $(t,x) \in [\tau,T] \times \supp(\mu^*(\tau)) \mapsto \F^{\omega,\tau}_t(x)$ is the unique solution of the Cauchy problem
\begin{equation}
\label{eq:Nonlocal_ODE_flow}
\left\{
\begin{aligned}
& \partial_t \F^{\omega,\tau}_t(x) = \INTDom{\BGamma^{\circ}_{ \left( t,\Phi^{v,u^*}_{(\tau,t)}[\mu^*(\tau)](x) \right)} \left( \Phi^{v,u^*}_{(\tau,t)}[\mu^*(\tau)](y) \right) \cdot \F^{\omega,\tau}_t(y) }{\R^d}{\mu^*(\tau)(y)} \\
& \hspace{1.75cm} + \left[ \D_x u^* \left( t , \Phi^{v,u^*}_{(\tau,t)}[\mu^*(\tau)](x) \right) + \D_x v [\mu^*(t)] \left( t , \Phi^{v,u^*}_{(\tau,t)}[\mu^*(\tau)](x) \right) \right] \cdot \F^{\omega,\tau}_t(x) \\
& \F^{\omega,\tau}_{\tau}(x) = \omega(x) - u^*(\tau,x). 
\end{aligned}
\right.
\end{equation}
\end{lem}

\begin{proof}
We start by computing the measures $\tilde{\mu}_{\tau}(\epsilon)$ as a function of $\mu^*(\tau)$ for all $\epsilon \in [0,\bar{\epsilon})$. By definition of the needle-like variation, it holds
\begin{equation*}
\tilde{\mu}_{\tau}(\epsilon) = \Phi^{v,\omega}_{(\tau-\epsilon,\tau)}[\mu^*(\tau-\epsilon)] \circ \Phi^{v,u^*}_{(\tau,\tau-\epsilon)}[\mu^*(\tau)](\cdot)_{\#} \mu^*(\tau)
\end{equation*}

Using Lebesgue's Differentiation Theorem, we obtain the following expansions at the first order with respect to $\epsilon$ 
\begin{equation*}
\begin{aligned}
\Phi^{v,\omega}_{(\tau-\epsilon,\tau)}[\mu^*(\tau-\epsilon)](y) & = y + \INTSeg{\left[ v[\tilde{\mu}_t(\epsilon)] \left( t,\Phi^{v,\omega}_{(\tau-\epsilon,t)}[\mu^*(\tau-\epsilon)](y) \right) + \omega \left( \Phi^{v,\omega}_{(\tau-\epsilon,t)}[\mu^*(\tau-\epsilon)](y) \right) \right]}{t}{\tau-\epsilon}{\tau}, \\
& = y + \epsilon \left( v[\mu^*(\tau)](\tau,y) + \omega(y) \right) + o(\epsilon), 
\end{aligned}
\end{equation*}
as well as 
\begin{equation*}
\begin{aligned}
\Phi^{v,u^*}_{(\tau,\tau-\epsilon)}[\mu^*(\tau)](y) & = y - \INTSeg{\left[ v[\mu^*(t)] \left( t,\Phi^{v,u^*}_{(t,\tau-\epsilon)}[\mu^*(\tau)](y) \right) + u^*\left( t , \Phi^{v,u^*}_{(t,\tau-\epsilon)}[\mu^*(\tau)](y) \right) \right]}{t}{\tau-\epsilon}{\tau}, \\
& = y - \epsilon \left( v[\mu^*(\tau)](\tau,y) + u^*(\tau,y) \right) + o(\epsilon). 
\end{aligned}
\end{equation*}

Chaining these two expressions together and recalling that $\omega(\cdot)$ and $v[\mu^*(\tau)](\tau,\cdot)$ are $C^1$-smooth, it holds
\begin{equation*}
\Phi^{v,\omega}_{(\tau-\epsilon,\tau)}[\mu^*(\tau-\epsilon)] \circ \Phi^{v,u^*}_{(\tau,\tau-\epsilon)}[\mu^*(\tau)](y) = y + \epsilon \left[ \omega(y) - u^*(\tau,y) \right] + o(\epsilon)
\end{equation*}
and we deduce the expression that will prove useful in the sequel
\begin{equation}
\label{eq:Perturbed_measure_tau}
\tilde{\mu}_{\tau}(\epsilon) = \left( \Id + \epsilon \left[ \omega(\cdot) - u^*(\tau,\cdot) \right] + o(\epsilon) \right)_{\#} \mu^*(\tau).
\end{equation}

We now want to obtain a similar expression but at some time $t \in [\tau,T]$. First, recall that $\tilde{\mu}_t(\epsilon) = \G^{\omega,\tau}_t(\epsilon,\cdot)_{\#} \mu^*(t)$ where 
\begin{equation}
\begin{aligned}
\G^{\omega,\tau}_t(\epsilon,\cdot) : x \in \supp(\mu^*(t)) \mapsto & \Phi^{v,u^*}_{(\tau,t)}[\tilde{\mu}_{\tau}(\epsilon)] \, \circ \, \Phi^{v,\omega}_{(\tau-\epsilon,\tau)}[\mu^*(\tau-\epsilon)] \\
\circ \, & \Phi^{v,u^*}_{(\tau,\tau-\epsilon)}[\mu^*(\tau)] \, \circ \, \Phi^{v,u^*}_{(t,\tau)}[\mu^*(t)](x).
\end{aligned}
\end{equation}
By \eqref{eq:Perturbed_measure_tau}, one has the following expansion
\begin{equation}\label{e-Oeps}
\begin{aligned}
& \hspace{0.5cm} \Phi^{v,u^*}_{(\tau,t)}[\tilde{\mu}_{\tau}(\epsilon)] \circ \Phi^{v,\omega}_{(\tau-\epsilon,\tau)}[\mu^*(\tau-\epsilon)] \circ \Phi^{v,u^*}_{(\tau,\tau-\epsilon)}[\mu^*(\tau)] \circ \Phi^{v,u^*}_{(t,\tau)}[\mu^*(t)](\cdot)  \\
& = \Phi^{v,u^*}_{(\tau,t)}[\tilde{\mu}_{\tau}(\epsilon)] \circ \left( \Phi^{v,u^*}_{(t,\tau)}[\mu^*(\tau)](\cdot) + \epsilon \left[ \omega \left(\Phi^{v,u^*}_{(t,\tau)}[\mu^*(\tau)](\cdot) \right) - u^* \left( \tau,\Phi^{v,u^*}_{(t,\tau)}[\mu^*(\tau)](\cdot) \right) \right] + o(\epsilon) \right) \\
& = \Phi^{v,u^*}_{(\tau,t)}[\tilde{\mu}_{\tau}(\epsilon)]\left( \Phi^{v,u^*}_{(t,\tau)}[\mu^*(\tau)](\cdot) \right) \\
& \hspace{0.45cm} + \epsilon \D_x \Phi^{v,u^*}_{(\tau,t)}[\mu^*(\tau)] \left( \Phi^{v,u^*}_{(t,\tau)}[\mu^*(\tau)](\cdot) \right) \left[ \omega \left(\Phi^{v,u^*}_{(t,\tau)}[\mu^*(\tau)](\cdot) \right) - u^* \left( \tau,\Phi^{v,u^*}_{(t,\tau)}[\mu^*(\tau)](\cdot) \right) \right] + o(\epsilon).
\end{aligned}
\end{equation}
since $\tilde{\mu}_{\tau}(\epsilon) \overset{W_1 \,}{\longrightarrow} \mu^*(\tau)$ as $\epsilon \downarrow 0$, $\mu \in \Pcal_c(\R^d) \mapsto \D_x \Phi_{(0,t)}^{v,u^*}[\mu](x)$ is continuous by \textbf{(F)}, and we are only interested in a Taylor expansion at the first order in $\epsilon$.

It then remains to compute the $O(\epsilon)$ term arising from $\Phi^{v,u^*}_{(\tau,t)}[\tilde{\mu}_{\tau}(\epsilon)](\cdot)$ in \eqref{e-Oeps}. Due to Proposition \ref{prop:Directional_derivative_flow}, the derivative of the non-local flow along directions induced by Lipschitz families of continuous and bounded maps exists. Recalling \eqref{eq:Perturbed_measure_tau}, this translates into
\begin{equation*}
\Phi^{v,u^*}_{(\tau,t)}[\tilde{\mu}_{\tau}(\epsilon)](y) = \Phi^{v,u^*}_{(\tau,t)}[\mu^*(\tau)](y) + \epsilon w_{\Phi}^{\omega,\tau}(t,y) + o(\epsilon), 
\end{equation*}
where $w_{\Phi}^{\omega,\tau}(t,y)$ is defined through \eqref{eq:Directional_derivative_charac} in the case where the non-local velocity field is given by $(t,x) \mapsto v[\mu^*(t)](t,x) + u^*(t,x)$.

Thus, we proved the pointwise Taylor expansion at the first order with $\epsilon$ :
\begin{equation*}
\G^{\omega,\tau}_t(\epsilon,x) = x + \epsilon \F^{\omega,\tau}_t \circ \Phi^{u^*,v}_{(t,\tau)}[\mu^*(t)](x) + o(\epsilon) 
\end{equation*}
for $\mu^*(t)$-almost every $x \in \R^d$, where 
\begin{equation}
\label{eq:FirstOrder_Needle_Expression_non-local}
\F^{\omega,\tau}_t : x \in \supp(\mu^*(\tau)) \mapsto \D_x \Phi^{v,u^*}_{(\tau,t)}[\mu^*(\tau)](x) \cdot \left[ \omega(x) - u^*(\tau,x) \right] + w_{\Phi}^{\omega,\tau}(t,x)
\end{equation}
is a continuous mapping for all $t \in [\tau,T]$.

A standard application of Proposition \ref{prop:Convergence_results}-$(iii)$ shows that this expansion holds in $L^2(\R^d,\R^d;\mu^*(t))$. One can then extend $\G^{\omega,\tau}_t(\cdot,\cdot)$ to $(-\bar{\epsilon},\bar{\epsilon})$ while preserving this expansion around $\epsilon = 0$ by defining e.g.
\begin{equation*}
\G^{\omega,\tau}_t(\epsilon,\cdot) \equiv \Id + \epsilon \F^{\omega,\tau}_t \circ \Phi^{v,u^*}_{(\tau,t)}[\mu^*(\tau)](\cdot) ~~ \text{for $\epsilon \in (-\bar{\epsilon},0]$}.
\end{equation*}
The existence of a constant $R_T^{\Phi}$ depending on $R_T,L_U$ and $v[\cdot](\cdot,\cdot)$ such that $\supp(\G^{\omega,\tau}_t(\epsilon,\cdot)_{\#}\mu^*(t)) \subset \overline{B(0,R_T^{\Phi})}$ follows from hypotheses \textbf{(F)} and \textbf{(B)}, which ensure the continuity and boundedness of the perturbation functions over $\overline{B(0,\tilde{R}_T)}$.

We finally prove the counterpart of Lemma \ref{lem:Differential_F} providing an ODE-type characterization for the perturbation induced by the needle-like variation in the non-local case. Recalling the definition of $(t,x) \mapsto \F^{\omega,\tau}_t(x)$ given in \eqref{eq:FirstOrder_Needle_Expression_non-local} and summing the ODE-type characterization of $t \mapsto w_{\Phi}^{\omega,\tau}(t,\cdot)$ and $\D_x \Phi^{v,u^*}_{(\tau,\cdot)}[\mu^*(\tau)](\cdot) \cdot \left[ \omega(\cdot) - u^*(\tau,\cdot) \right]$, we recover
\begin{equation*}
\left\{
\begin{aligned}
& \partial_t \F^{\omega,\tau}_t(x) = \INTDom{\BGamma^{\circ}_{ \left( t,\Phi^{v,u^*}_{(\tau,t)}[\mu^*(\tau)](x) \right)} \left( \Phi^{v,u^*}_{(\tau,t)}[\mu^*(\tau)](y) \right) \cdot \F^{\omega,\tau}_t(y) }{\R^d}{\mu^*(\tau)(y)} \\
& \hspace{1.75cm} + \left[ \D_x u^* \left( t , \Phi^{v,u^*}_{(\tau,t)}[\mu^*(\tau)](x) \right) + \D_x v [\mu^*(t)] \left( t , \Phi^{v,u^*}_{(\tau,t)}[\mu^*(\tau)](x) \right) \right] \cdot \F^{\omega,\tau}_t(x)  \\
& \F^{\omega,\tau}_{\tau}(x) = \omega(x) - u^*(\tau,x), 
\end{aligned}
\right.
\end{equation*}
which concludes the proof of our result.
\end{proof}

In the development of Steps 2, 3 and 4, we do not need to take into account the explicit dependence of the flows with respect to their starting measures. We shall henceforth write $\Phi^{v,u^*}_{(\cdot,\cdot)}(\cdot) \equiv \Phi^{v,u^*}_{(\cdot,\cdot)}[\mu^*(\cdot)](\cdot)$ for clarity and conciseness.

\begin{flushleft}
\underline{\textbf{Step 2 : First-order optimality condition}}
\end{flushleft}

In the framework of Problem $(\pazocal{P})$, the optimality of $u^*(\cdot)$ writes
\begin{equation}
\label{eq:Optimality_General}
\begin{aligned}
\frac{\varphi(\tilde{\mu}_T(\epsilon)) - \varphi(\mu^*(T))}{\epsilon} & + \frac{1}{\epsilon} \INTSeg{ \left[ L(\tilde{\mu}_t(\epsilon),\omega) - L(\mu^*(t),u^*(t)) \right]}{t}{\tau - \epsilon}{\tau} \\
& + \frac{1}{\epsilon} \INTSeg{\left[ L(\tilde{\mu}_t(\epsilon),u^*(t)) - L(\mu^*(t),u^*(t)) \right] }{t}{\tau}{T} ~ \geq 0,
\end{aligned}
\end{equation}
for all $\epsilon \in [0,\bar{\epsilon})$. 

The first order perturbation corresponding to the final cost $\varphi(\cdot)$ has already been treated in \eqref{eq:Optimality_endcost}-\eqref{eq:Gradient_terminal_cost}, Section \ref{subsection:SimplePMP}. We study the integral terms arising from the running cost. Remark first that it holds
\begin{equation*}
\lim\limits_{\epsilon \downarrow 0} \left[ \frac{1}{\epsilon} \INTSeg{ \left[ L(\tilde{\mu}_t(\epsilon),\omega) - L(\mu^*(t),u^*(t)) \right]}{t}{\tau - \epsilon}{\tau} \right] = L(\mu^*(\tau),\omega) - L(\mu^*(\tau),u^*(\tau)),
\end{equation*}
by the Lebesgue Differentiation Theorem, since $\tilde{\mu}_t(\epsilon) \overset{W_1 \,}{\longrightarrow} \mu^*(t)$ as $\epsilon \downarrow 0$ for all $t \in [0,T]$ and since $\tau$ is a Lebesgue point of $u^*(\cdot)$.

Equivalently to the proof of the PMP for Problem $(\pazocal{P}_1)$, the perturbed measures are uniformly supported in a compact set. Thus, under hypotheses \textbf{(L)} and recalling that the function $\epsilon \in (-\bar{\epsilon},\bar{\epsilon}) \mapsto L(\tilde{\mu}_t(\epsilon),u^*(t))$ is differentiable at $\epsilon = 0$ for $\Lcal^1$-almost every $t \in [\tau,T]$ by hypothesis \textbf{(D)}, the chain rule of Proposition \ref{prop:Chainrule} can be applied to the running cost to obtain 
\begin{equation} 
\label{eq:Chainrule_runningcost}
\lim\limits_{\epsilon \downarrow 0} \left[ \frac{1}{\epsilon} \left[ L(\tilde{\mu}_t(\epsilon),u^*(t)) - L(\mu^*(t),u^*(t)) \right] \right] = \INTDom{\langle \bar{\Bgamma}^{\circ}_L(t,x) , \F^{\omega,\tau}_t \circ \Phi^{v,u^*}_{(t,\tau)}(x) \rangle}{\R^d}{\mu^*(t)(x)},
\end{equation}
where $\bar{\Bgamma}^{\circ}_L(t,\cdot) \in L^2(\R^d,\R^d;\mu^*(t))$ is the barycenter of $\Bpartial^{\circ}_{\mu} L(\mu^*(t),u^*(t))$ for $\Lcal^1$-almost every $t \in [0,T]$.

Moreover, the uniform compactness of the supports of the perturbed measures and hypothesis \textbf{(L)} imply that the left hand side in \eqref{eq:Chainrule_runningcost} is uniformly bounded by a function in $L^1([0,T],\R_+)$ for any $\epsilon \in (0,\bar{\epsilon})$. Therefore, it holds by an application of Lebesgue Dominated Convergence Theorem that
\begin{equation*}
\lim\limits_{\epsilon \downarrow 0} \left[ \frac{1}{\epsilon} \INTSeg{\left[ L(\tilde{\mu}_t(\epsilon),u^*(t)) - L(\mu^*(t),u^*(t)) \right]}{t}{\tau}{T} \right] =\INTSeg{\INTDom{\langle \bar{\Bgamma}^{\circ}_L(t,x) , \F^{\omega,\tau}_t \circ \Phi^{v,u^*}_{(t,\tau)}(x) \rangle}{\R^d}{\mu^*(t)(x)}}{t}{\tau}{T}.
\end{equation*}
Thus, the optimality of $(u^*(\cdot),\mu^*(\cdot))$ translates into the first-order condition
\begin{equation} 
\label{eq:First_Order_General}
\begin{aligned}
\INTDom{\langle \bar{\Bgamma}^{\circ}_{\varphi}(x) , \F^{\omega,\tau}_T \circ \Phi^{v,u^*}_{(T,\tau)}(x) \rangle}{\R^d}{\mu^*(T)(x)} & + \left[ L(\mu^*(\tau),\omega) - L(\mu^*(\tau),u^*(\tau)) \right] \\
& + \INTSeg{\INTDom{\langle \bar{\Bgamma}^{\circ}_L(t,x) , \F^{\omega,\tau}_t \circ \Phi^{v,u^*}_{(t,\tau)}(x) \rangle}{\R^d}{\mu^*(t)(x)}}{t}{\tau}{T} ~ \geq ~ 0,
\end{aligned}
\end{equation}
for any couple of needle parameters $(\omega,\tau)$. 

\begin{flushleft}
\textbf{\underline{Step 3 : Backward dynamics and Pontryagin maximization condition}}
\end{flushleft}

We now build a solution $\nu^*(\cdot) \in \Lip([0,T],\Pcal_c(\R^{2d}))$ to the system of continuity equations
\begin{equation}
\label{eq:GeneralPMP_HamiltonianDynamics}
\left\{ 
\begin{aligned}
& \partial_t \nu^*(t) + \nabla \cdot (\V^*[\nu^*(t)](t,\cdot,\cdot) \nu^*(t)) = 0 ~~ & \text{in $[0,T] \times \R^{2d}$}, \\
& \pi^1_{\#} \nu^*(t) = \mu^*(t) ~~ & \text{for all $t \in [0,T]$}, \\
& \nu^*(T)= (\Id \times (-\bar{\Bgamma}^{\circ}_{\varphi}))_{\#} \mu^*(T),
\end{aligned}
\right.
\end{equation}
associated to the non-local vector field 
\begin{equation*}
\begin{aligned}
\V^*[\nu^*(t)] : & (t,x,r) \in [0,T] \times \supp(\nu) \mapsto \\
& \hspace{-0.1cm} \begin{pmatrix} v[\pi^1_{\#} \nu^*(t)](t,x) + u^*(t,x) \\ \bar{\Bgamma}^{\circ}_L(t,x) - \BGamma^{\circ}_v[\nu^*(t)](t,x) - \D_x u^*(t,x)^{\top} r - \D_x v[\mu^*(t)](t,x)^{\top}r) \end{pmatrix}
\end{aligned}
\end{equation*}
where $\BGamma^{\circ}_v[\cdot](\cdot,\cdot)$ is given by 
\begin{equation}
\label{eq:BGamma_def}
\BGamma^{\circ}_v[\nu](t,x) = \INTDom{\left( \BGamma^{\circ}_{(t,y)}(x) \right)^{\top} p \,}{\R^{2d}}{\nu(y,p)}
\end{equation}
for any $(\nu,t,x) \in \Pcal_c(\R^d) \times [0,T] \times \supp(\pi^1_{\#}\nu)$, with $\BGamma^{\circ}_{(\cdot,\cdot)}(\cdot)$ defined as in Theorem \ref{thm:General_PMP}. 

As in Lemma \ref{lem:Wellposedness_Simple} of Section \ref{subsection:SimplePMP}, we build a solution $\nu^*(\cdot)$ of \eqref{eq:GeneralPMP_HamiltonianDynamics} by making use of the cascaded structure of the system. We then show that this solution is such that the map $\K_{\omega,\tau}(\cdot)$ defined in this context by
\begin{equation}
\label{eq:General_K}
\begin{aligned}
\K_{\omega,\tau} : t \in [\tau,T] \mapsto & \INTDom{\langle r , \F^{\omega,\tau}_t \circ \Phi^{v,u^*}_{(t,\tau)}(x) \rangle}{\R^{2d}}{\nu^*(t)(x,r)} + \left[ L(\mu^*(\tau),u^*(\tau)) - L(\mu^*(\tau),\omega) \right] \\
& - \INTSeg{\INTDom{\langle \bar{\Bgamma}^{\circ}_L(t,x) , \F^{\omega,\tau}_s \circ \Phi^{v,u^*}_{(s,\tau)}(x) \rangle}{\R^{2d}}{\mu^*(s)(x)}}{s}{\tau}{t} \\
\end{aligned}
\end{equation}
is constant over $[\tau,T]$. 

\begin{lem}[Well-posedness of solutions of \eqref{eq:GeneralPMP_HamiltonianDynamics}]
\label{lem:Wellposedness_General}
Let $(u^*(\cdot),\mu^*(\cdot))$ be an optimal pair control-trajectory for $(\pazocal{P})$. We consider the family of maps $(t,x,r) \in [0,T] \times \R^{2d} \mapsto \Psi_t(x,r) \in \R^d$, defined to be the solution of
\begin{equation}
\label{eq:Marginal_TransportEq_General}
\left\{
\begin{aligned}
& \partial_t w(t,x,r) = \bar{\Bgamma}^{\circ}_L(t,\Phi^{v,u^*}_{(T,t)}(x)) - \D_x u^*(t,\Phi^{v,u^*}_{(T,t)}(x))^{\top} w(t,x,r) - \D_x v[\mu^*(t)](t,\Phi^{v,u^*}_{(T,t)}(x))^{\top} w(t,x,r) \\
& \hspace{4.25cm} - \INTDom{\BGamma^{\circ}_{\left( t , \Phi^{v,u^*}_{(T,t)}(y) \right)} \left( \Phi^{v,u^*}_{(T,t)}(x) \right)^{\top} w(t,y,p) \, }{\R^{2d}}{(\Id \times (-\bar{\Bgamma}^{\circ}_{\varphi}))_{\#} \mu^*(T) (y,p)}, \\
& w(T,x,r) = r
\end{aligned}
\right.
\end{equation}
For $\mu^*(T)$-almost every $x \in \R^d$, we define the curves of measures $\sigma^*_x : t \mapsto \Psi_t(x,\cdot)_{\#} \delta_{(-\bar{\Bgamma}^{\circ}_{\varphi}(x))}$ and denote by $\V_x^*(\cdot,\cdot)$ the corresponding non-local vector fields describing their evolution. 

Then $\nu^* : t \mapsto (\Phi^{v,u^*}_{(T,t)} , \Id)_{\#} \nu^*_T(t)$ solves \eqref{eq:GeneralPMP_HamiltonianDynamics} with $\nu^*_T(t) = \int \sigma^*_x(t) \textnormal{d} \mu^*(T)(x)$. Moreover, there exists two constants $R'_T,L_T' > 0$ such that
\begin{equation*}
\supp(\nu^*(t)) \subset \overline{B_{2d}(0,R'_T)} ~~ \text{and} ~~ W_1(\nu^*(t),\nu^*(s)) \leq L_T' |t-s|,
\end{equation*}  
for all $s,t \in [0,T]$
\end{lem}

\begin{proof}
We denote by $\Omega$ the compact\footnote{Recall that $\bar{\Bgamma}^{\circ}_{\varphi}(\cdot)$ is a continuous map by hypothesis \textbf{(B)}.} subset $\supp((\Id \times (-\bar{\Bgamma}^{\circ}_{\varphi})_{\#}\mu^*(T))$ of $\R^{2d}$. We first show that \eqref{eq:Marginal_TransportEq_General} admits a unique continuous solution $(t,x,r) \in [0,T] \times \Omega \mapsto \Psi_t(x,r) \in \R^d$. This can be done by reproducing the strategy of the proof of Proposition \ref{prop:Directional_derivative_flow} which consists in defining a weighted $C^0([0,T] \times \Omega)$-norm for which the right-hand side of \eqref{eq:Marginal_TransportEq_General} is contracting and applying Banach's Fixed Point Theorem. Notice that here, the coupling between the non-local flows arising from the integral term in \eqref{eq:Marginal_TransportEq_General} requires us to use explicitly the continuity of the right-hand side with respect to $x$. In Lemma \ref{lem:Wellposedness_Simple}, all the backward Cauchy problems were independent and we did not need any regularity assumption on $x$ for the proof to work.

Since $[0,T] \times \Omega$ is compact, $(t,x,r) \mapsto \Psi_t(x,r)$ is bounded. This implies by \eqref{eq:Marginal_TransportEq_General} that $t \mapsto \Psi_t(x,r)$ is Lipschitz for all $(x,r) \in \Omega$. Moreover, a direct application of Gr\"onwall Lemma along with \eqref{eq:Marginal_TransportEq_General} allows to show that for all $(t,x) \in [0,T] \times \pi^1(\Omega)$, ones has
\begin{equation*}
|\Psi_t(x,r_2) - \Psi_t(x,r_1) | \leq C |r_2-r_1|
\end{equation*}
for all $(r_1,r_2) \in \pi^2(\Omega)$. Hence, we showed that $(t,r) \in [0,T] \times \pi^2(\Omega) \mapsto \Psi_t(x,r)$ is Lipschitz for all $x \in \pi^1(\Omega)$.

\begin{comment}
Now, we want to prove that these maps are uniformly bounded. To do so, we introduce the function
%
\begin{equation*}
\Wcal : t \in [0,T] \mapsto \max \left\{ \frac{1}{2}|\Psi_t(x,r)|^2 ~\text{s.t.}~ (x,r) \in \Omega \right\}. 
\end{equation*}
%
Since $(x,r) \mapsto \frac{1}{2}|\Psi_t(x,r)|^2$ is continuous for all  $t \in [0,T]$ and since $\Omega$ is compact, we can invoke Danskin's Theorem to write 
%
\begin{equation*}
\derv{}{t} \Wcal(t) = \max \left\{ \left\langle \partial_t \Psi_t(x_*,r_*) , \Psi_t(x_*,r_*) \right\rangle ~\text{s.t.}~ (x_*,r_*) \in \Omega_t \right\},
\end{equation*}
%
where $\Omega_t = \{ (x_*,r_*) \in \Omega ~\text{s.t.}~ \frac{1}{2}|\Psi_t(x_*,r_*)|^2 = \Wcal(t) \}$. Using the expression given in \eqref{eq:Marginal_TransportEq_General} for $\partial_t \Psi_t(x_*,r_*)$, Cauchy-Schwarz inequality and the fact that all the quantities involved are bounded over $[0,T] \times \Omega$, we recover that there exists two constants $C_1,C_2 > 0$ such that 
%
\begin{equation*}
\derv{}{t} \Wcal(t) \leq C_1 + C_2 \Wcal(t) ~~ \text{for all $t \in [0,T]$}.
\end{equation*}

This implies by Gr\"onwall Lemma that $\Wcal(t) \leq \Wcal(T)$ for all $t \in [0,T]$, which is bounded by definition of $\Omega$. This allows to conclude, up to performing another Gr\"onwall estimate, that the maps $(t,r) \mapsto \Psi_t(x,r)$ are Lipschitz over $[0,T] \times \pi^2(\Omega)$ uniformly with respect to $x \in \mu^*(T)$.
\end{comment}

Therefore, carrying out same computations as in Lemma \ref{lem:Wellposedness_Simple}, we show that the curves of measures $t \mapsto \sigma^*_x(t) = \Psi_t(x,\cdot)_{\#} \delta_{-\bar{\Bgamma}^{\circ}_{\varphi(x)}}$ are well-defined, uniformly compactly supported and Lipschitz in the $W_1$-metric for $\mu^*(T)$-almost every $x \in \R^d$. This implies the existence of two constants $R_T',L_T' > 0$ such that it holds
\begin{equation*}
\supp(\nu^*(t)) \subset \overline{B_{2d}(0,R'_T)} ~~ \text{and} ~~ W_1(\nu^*(t),\nu^*(s)) \leq L_T' |t-s| ~~ \text{for all $s,t \in [0,T]$}.
\end{equation*}  

Moreover, for any $\xi \in C^{\infty}_c(\R^{2d},\R)$, it holds that
\begin{equation*}
\begin{aligned}
& \derv{}{t} \left[ \INTDom{\xi(x,r)}{\R^{2d}}{\nu^*(t)(x,r)} \right] =  \derv{}{t} \left[ \INTDom{\INTDom{\xi(\Phi^{v,u^*}_{(T,t)}(x),r)}{\R^d}{\sigma_x^*(t)(r)}}{\R^d}{\mu^*(T)(x)} \right] \\
= & \INTDom{\INTDom{\langle \nabla_x \xi(\Phi^{v,u^*}_{(T,t)}(x),r) , v[\mu^*(t)](t,\Phi^{v,u^*}_{(T,t)}(x)) \rangle}{\R^d}{\sigma^*_x(t)(r)}}{\R^d}{\mu^*(T)(x)} \\
+ & \INTDom{\INTDom{\langle \nabla_x \xi(\Phi^{v,u^*}_{(T,t)}(x),r) , u^*(t,\Phi^{v,u^*}_{(T,t)}(x))\rangle}{\R^d}{\sigma^*_x(t)(r)}}{\R^d}{\mu^*(T)(x)} \\
+ & \INTDom{\INTDom{\langle \nabla_r \xi(\Phi^{v,u^*}_{(T,t)}(x),r) , \bar{\Bgamma}^{\circ}_L(t,\Phi^{v,u^*}_{(T,t)}(x)) \rangle}{\R^d}{\sigma^*_x(t)(r)}}{\R^d}{\mu^*(T)(x)} \\
- & \INTDom{\INTDom{\langle \nabla_r \xi(\Phi^{v,u^*}_{(T,t)}(x),r) , \D_x u^*(t,\Phi^{v,u^*}_{(T,t)}(x))^{\top} r \rangle}{\R^d}{\sigma^*_x(t)(r)}}{\R^d}{\mu^*(T)(x)} \\
- & \INTDom{\INTDom{\langle \nabla_r \xi(\Phi^{v,u^*}_{(T,t)}(x),r) , \D_x v[\mu^*(t)](t,\Phi^{v,u^*}_{(T,t)}(x))^{\top}r \rangle}{\R^d}{\sigma^*_x(t)(r)}}{\R^d}{\mu^*(T)(x)} \\
- & \INTDom{\INTDom{ \left\langle \nabla_r \xi(\Phi^{v,u^*}_{(T,t)}(x),r) , \BGamma^{\circ}_v[\nu^*(t)](t,\Phi^{v,u^*}_{(T,t)}(x)) \right\rangle}{\R^d}{\sigma^*_x(t)(r)}}{\R^d}{\mu^*(T)(x)} \\
= & \INTDom{\langle \nabla_x \xi(x,r) , v[\mu^*(t)](t,x) + u^*(t,x)\rangle}{\R^{2d}}{\nu^*(t)(x,r)} \\
+ & \INTDom{\langle \nabla_r \xi(x,r) , \bar{\Bgamma}^{\circ}_L(t,x) - \D_x u^*(t,x)^{\top} r  - \D_x v[\mu^*(t)](t,x)^{\top} r \rangle}{\R^{2d}}{\nu^*(t)(x,r)} \\
- & \INTDom{\langle \nabla_r \xi(x,r) , \BGamma^{\circ}_v[\nu^*(t)](t,x) \rangle}{\R^{2d}}{\nu^*(t)(x,r)}
\end{aligned}
\end{equation*}
where we used the fact that
\begin{equation*}
\begin{aligned}
\BGamma^{\circ}_v[\nu^*(t)](t,\Phi^{v,u^*}_{(T,t)}(x)) = \INTDom{\BGamma^{\circ}_{\left( t , \Phi^{v,u^*}_{(T,t)}(y) \right)} \left( \Phi^{v,u^*}_{(T,t)}(x) \right)^{\top} p \, }{\R^{2d}}{\nu_T^*(t)(y,p)}&  \\
= \INTDom{\INTDom{\BGamma^{\circ}_{\left( t , \Phi^{v,u^*}_{(T,t)}(y) \right)} \left( \Phi^{v,u^*}_{(T,t)}(x) \right)^{\top} p \, }{\R^d}{\sigma^*_y(t)(p)}}{\R^d}{\mu^*(T)(y)}&  \\
= \INTDom{\BGamma^{\circ}_{\left( t , \Phi^{v,u^*}_{(T,t)}(y) \right)} \left( \Phi^{v,u^*}_{(T,t)}(x) \right)^{\top} \Psi_t(y,p) \, }{\R^{2d}}{((\Id \times (-\bar{\Bgamma}^{\circ}_{\varphi})_{\#}\mu^*(T))(y,p)} & 
\end{aligned}
\end{equation*}
We therefore recover that $t \mapsto \nu^*(t)$ solves \eqref{eq:GeneralPMP_HamiltonianDynamics}, which ends the proof.
\end{proof}

\begin{lem}
\label{lem:Backward_VectorField_General}
The map $t \mapsto \K_{\omega,\tau}(t)$ defined in \eqref{eq:General_K} is constant over $[\tau,T]$ for any couple of needle parameters $(\omega,\tau)$. 
\end{lem}

\begin{proof}
This proof follows the same steps as in the proof of Lemma \ref{lem:Backward_VectorField_Simple}, the  difference lying in the fact that the flows $(\Phi^{v,u^*}_{(0,t)}(\cdot))_{t \in [0,T]}$ are associated to the non-local PDE. It can be verified again as in the proof of Lemma \ref{lem:Backward_VectorField_Simple} that $t \mapsto \K_{\omega,\tau}(t)$ is Lipschitz. We compute $\derv{}{t} \K_{\omega,\tau}$ using \eqref{eq:NonlocalPDE_distributions2} as in \eqref{eq:Time_derivative_K} while plugging in the expressions for $\F^{\omega,\tau}_t(\cdot)$ and its time-derivative provided by Lemma \ref{lem:Perturbed_state_non-local}.

\begin{equation*}
\begin{aligned}
\derv{}{t} \K_{\omega,\tau}(t) & = \INTDom{ \INTDom{\left\langle r , \partial_t \F^{\omega,\tau}_t \circ \Phi^{v,u^*}_{(T,\tau)}(x) \right\rangle}{\R^d}{\sigma^*_x(t)(r)} }{\R^d}{\mu^*(T)(x)}  \\
& \hspace{0.45cm} + \INTDom{ \INTDom{\left\langle \V_x^*(t,r) , \F^{\omega,\tau}_t \circ \Phi^{v,u^*}_{(T,\tau)}(x) \right\rangle}{\R^d}{\sigma^*_x(t)(r)} }{\R^d}{\mu^*(T)(x)} \\
& \hspace{0.45cm} - \INTDom{\langle \bar{\Bgamma}^{\circ}_L ( t, \Phi^{v,u^*}_{(T,t)}(x) ) , \F^{\omega,\tau}_t \circ \Phi^{v,u^*}_{(T,\tau)}(x) \rangle}{\R^{2d}}{\mu^*(T)(x)} \\
& = \INTDom{\left\langle r , \INTDom{ \BGamma^{\circ}_{\left( t , \Phi^{v,u^*}_{(T,t)}(x) \right)} \left( \Phi^{v,u^*}_{(T,t)}(y) \right) \F^{\omega,\tau}_t \circ \Phi^{v,u^*}_{(T,\tau)}(x) }{\R^d}{\mu^*(T)(y)} \right\rangle}{\R^{2d}}{\nu^*_T(t)(x,r)} \\ 
& \hspace{0.45cm} - \INTDom{\left\langle \INTDom{ \BGamma^{\circ}_{\left( t , \Phi^{v,u^*}_{(T,t)}(y) \right)} \left( \Phi^{v,u^*}_{(T,t)}(x) \right)^{\top} p \, }{\R^{2d}}{\nu^*_T(t)(y,p)} , \F^{\omega,\tau}_t \circ \Phi^{v,u^*}_{(T,\tau)}(x) \right\rangle}{\R^{2d}}{\nu^*_T(t)(x,r)} \\
& = 0, 
\end{aligned}
\end{equation*}
by plugging in the expressions of $\partial_t \F^{\omega,\tau}_t(\cdot)$ and $\V^*_x(t,\cdot)$. The two quantities are shown to be equal due to the uniform boundedness of the integrands given by \textbf{(B)} and Fubini-Tonelli theorem. This altogether leads to $\derv{}{t} \K_{\omega,\tau}(t) = 0$ for $\Lcal^1$-almost every $t \in [\tau,T]$ and thus to $\K_{\omega,\tau}(\cdot)$ being constant over $[\tau,T]$.
\end{proof}

\begin{flushleft}
\textbf{\underline{Step 4 : Proof of the Pontryagin Maximum Principle for $(\pazocal{P})$}}
\end{flushleft}

We proved in Lemma \ref{lem:Wellposedness_General} that there exists a curve $\nu^* \in \Lip([0,T],\Pcal_c(\R^{2d}))$ solution of \eqref{eq:GeneralPMP_HamiltonianDynamics} along with a constant $R \equiv R_T' > 0$ such that $\supp(\nu^*(\cdot)) \subset \overline{B_{2d}(0,R)}$. The non-local velocity field $\V^*[\nu^*(\cdot)](\cdot,\cdot,\cdot)$ is defined for $\Lcal^1 \times \nu^*(\cdot)$-almost every $(t,x,r) \in [0,T] \times \overline{B_{2d}(0,R)}$ by
\begin{equation*}
\V^*[\nu^*(t)](t,x,r) = \begin{pmatrix}
\bar{\Bgamma}^{\circ}_L(t,x) - \D_x u^*(t,x)^{\top} r - \D_x v[\pi^1_{\#}\nu^*(t)](t,x)^{\top} r - \BGamma^{\circ}_v [\nu^*(t)](t,x) \\
v[\pi^1_{\#} \nu^*(t)](t,x) + u^*(t,x)
\end{pmatrix}.
\end{equation*}

We define the infinite dimensional Hamiltonian $\H(\cdot,\cdot,\cdot)$ of the system by
\begin{equation*}
\H : (t,\nu,\omega) \in [0,T] \times \Pcal_c(\R^{2d}) \times U \mapsto \INTDom{\langle r , \omega(x) \rangle}{\R^{2d}}{\nu(x,r)} - L(\pi^1_{\#} \nu , \omega).
\end{equation*}
along with its compactification $\H_c(\cdot,\cdot,\cdot)$ given by \eqref{eq:General_Hamiltonian_Compact}. 

Furthermore, we proved in Lemma \ref{lem:Backward_VectorField_General} that the solution $\nu^*(\cdot)$ that we built is such that the map $\K_{\omega,\tau}(\cdot)$ defined in \eqref{eq:General_K} is constant over $[\tau,T]$ for any couple of needle parameters $(\omega,\tau)$. Hence, it holds in particular that $\K_{\omega,\tau}(\tau) = \K_{\omega,\tau}(T)$ which is a non-positive quantity by the first-order optimality condition \eqref{eq:First_Order_General}. This fact implies that
\begin{equation*}
\INTDom{\langle r , \omega(x) - u^*(\tau,x) \rangle}{\R^{2d}}{\nu^*(\tau)(x,r)} - \left[ L(\mu^*(t),\omega) - L(\mu^*(t),u^*(t)) \right] \leq 0,
\end{equation*}
for all $\omega \in U$ and $\tau \in [0,T]$ Lebesgue point of $v[\mu^*(\cdot)](\cdot,\cdot)+u^*(\cdot)$. This inequality rewrites as the Pontryagin Maximization condition \eqref{eq:General_Maximization} : 
\begin{equation*}
\H_c(t,\nu^*(t),u^*(t)) = \max\limits_{\omega \in U} \left[ \H_c(t,\nu^*(t),\omega\right]
\end{equation*}
for $\Lcal^1$-almost every $t \in [0,T]$.

Finally, one recognizes the pseudo-Hamiltonian structure $\V^*[\nu^*(t)](t,x,r) = \J_{2d} \tilde{\nabla}_{\nu} \H_c(t,\nu^*(t),u^*(t))$ for $\Lcal^1 \times \nu^*(\cdot)$-almost every $(t,x,r) \in [0,T] \times \overline{B_{2d}(0,R)}$ where the map $\tilde{\nabla}_{\nu} \H_c (t,\nu^*(t),u^*(t))(\cdot,\cdot)$ is precisely the non-local velocity field $\V^*[\nu^*](t,\cdot,\cdot)$ for $\Lcal^1$-almost every $t \in [0,T]$. This concludes our proof of the Pontryagin Maximum Principle for $(\pazocal{P})$.

%%%%%%%%%%%%%%%%%%%%%%%%%%%%%%%%%%%%%%%%%%%%%%%%%%%%%%%%%%%%%%%%%%%%%%%%%%%%%%%%%%
%                               NEW SECTION AHEAD                                %
%%%%%%%%%%%%%%%%%%%%%%%%%%%%%%%%%%%%%%%%%%%%%%%%%%%%%%%%%%%%%%%%%%%%%%%%%%%%%%%%%%

\section{Examples}
\label{section:Examples}

The aim of the general result stated in Theorem \ref{thm:General_PMP} is to provide first-order necessary optimality conditions that are adapted to a wide range of functionals. We give in the following Propositions some examples of classical functionals that are encompassed in hypotheses \textbf{(H)} and compute the minimal selection in their Wasserstein subdifferential.

\begin{prop}[Subdifferential of a smooth integral functional]
\label{prop:Subdifferential_Integral}
\hfill \\ Let $V \in C^1(\R^d,\R)$ and $K \subset \R^d$ be a compact set. Define $\Vcal : \mu \in \Pcal(K) \mapsto \INTDom{V(x)}{\R^d}{\mu(x)}$. Then the functional $\Vcal(\cdot)$ is regular at any $\mu \in \Pcal(K)$ in the sense of Definition \ref{def:Regular}, Lipschitz in the $W_1$-metric. Moreover, the minimal selection $\Bpartial^{\circ} \Vcal(\mu)$ in its extended subdifferential at $\mu$ is a \textnormal{classical strong subdifferential} induced by a map and given explicitly by $\Bpartial^{\circ} \Vcal(\mu) = (\Id \times \nabla V)_{\#}\mu$. 
\end{prop}

\begin{proof}
See e.g. \cite[Proposition 10.4.2]{AGS}.
\end{proof}

\begin{rmk}
Taking any power $\alpha > 0$ of $\Vcal(\cdot)$ yields the same results provided that the functional $x \mapsto x^{\alpha}$ is differentiable at $\Vcal(\mu)$. In which case, the minimal selection in the extended subdifferential is induced by the map
\begin{equation}
\nabla_{\mu} (\Vcal^{\alpha})(\mu) : x \in \supp(\mu) \mapsto \alpha \Vcal(\mu)^{\alpha-1} \nabla V(x).
\end{equation}
\end{rmk}

\begin{prop}[Subdifferential of the variance functional]
\label{prop:Subdifferential_Variance}
Let $K \subset \R^d$ be a compact set and define the \textnormal{variance functional} by 
\begin{equation*}
\Var : \mu \in \Pcal(K) \mapsto \frac{1}{2} \INTDom{|x-\bar{\mu}|^2}{\R^d}{\mu(x)} = \frac{1}{2} \INTDom{|x|^2}{\R^d}{\mu(x)} - \frac{1}{2} |\bar{\mu}|^2
\end{equation*}
where $\bar{\mu} = \INTDom{x \,}{\R^d}{\mu(x)}$ denotes the average of the measure $\mu$.

Then, the functional $\Var(\cdot)$ is regular at any $\mu \in \Pcal(K)$, Lipschitz in the $W_1$-metric and the minimal selection $\Bpartial^{\circ} \Var(\mu)$ in its extended subdifferential is a  \textnormal{classical strong subdifferential} induced by the map $\nabla_{\mu} \Var(\mu) : x \in \supp(\mu) \mapsto x - \bar{\mu}$. 
\end{prop}

\begin{proof}
It is clear by definition of the variance functional that it is bounded from below over $\Pcal(K)$. Moreover, an application of the \textit{Kantorovich-Rubinstein} duality formula \eqref{eq:Kantorovich_duality} yields the Lipschitzianity in the $W_1$-metric. The regularity in the sense of Definition \ref{def:Regular} is a consequence of the convexity along Wasserstein geodesics of the variance functional (see \cite[Lemma 10.3.8]{AGS}).

We now show that $x \mapsto x - \bar{\mu}$ is in the classical strong subdifferential of the variance functional at $\mu \in \Pcal(K)$. For any $\nu \in \Pcal(K)$ and $\Bmu \in \Gamma(\mu,\nu)$, it holds that
\begin{equation*}
\begin{aligned}
& \INTDom{\langle x_1 - \bar{\mu} , x_2 - x_1 \rangle}{\R^{2d}}{\Bmu(x_1,x_2)} \\
= & \INTDom{\langle x_1 , x_2 \rangle}{\R^{2d}}{\Bmu(x_1,x_2)} - \INTDom{|x_1|^2}{\R^d}{\mu(x_1)} + |\bar{\mu}|^2 - \langle \bar{\mu} , \bar{\nu} \rangle, \\
\leq & \frac{1}{2} \INTDom{|x_2|^2}{\R^d}{\nu(x_2)} - \frac{1}{2} \INTDom{|x_1|^2}{\R^d}{\mu(x_1)} + |\bar{\mu}|^2 - \langle \bar{\mu} , \bar{\nu} \rangle , \\
\leq & \Var(\nu) - \Var(\mu) + \frac{1}{2} |\bar{\mu} - \bar{\nu}|^2.
\end{aligned}
\end{equation*}

Moreover, one can estimate the quantity $|\bar{\mu} - \bar{\nu}|^2$ as follows:
\begin{equation*} 
\begin{aligned}
|\bar{\mu} - \bar{\nu}|^2 & \leq \left( \INTDom{|x_1-x_2|}{\R^{2d}}{\Bmu(x_1,x_2)} \right)^2 \\ & \leq ~ \INTDom{|x_1-x_2|^2}{\R^{2d}}{\Bmu(x_1,x_2)} = W_{2,\Bmu}^2(\mu,\nu) = o(W_{2,\Bmu}(\mu,\nu)),
\end{aligned}
\end{equation*} 
since $\Bmu \in \Gamma(\mu,\nu)$, invoking Jensen's Inequality and the definition of $W_{2,\Bmu}(\cdot,\cdot)$ given in Definition \ref{def:Subdifferentials}.

Therefore, we conclude that for any $\nu \in \Pcal(K)$ and any $\Bmu \in \Gamma(\mu,\nu)$ it holds 
\begin{equation*}
\Var(\nu) - \Var(\mu) \geq \INTDom{\langle x_1 - \bar{\mu} , x_2 - x_1 \rangle}{\R^{2d}}{\Bmu(x_1,x_2)} + o(W_{2,\Bmu}(\mu,\nu)). 
\end{equation*}
which is equivalent to $x \in \supp(\mu) \mapsto x - \bar{\mu}$ being a classical strong subdifferential at $\mu$.

Now, take in particular $\nu \equiv \nu_s = (\Id + s \xi)_{\#} \mu$ for some small $s > 0$ and $\xi \in C^{\infty}_c(\R^d)$ such that $\supp(\nu_s) \subset \Pcal(K)$. It then holds
\begin{equation*}
+\infty > \lim\limits_{s \downarrow 0} \left[ \frac{\Var((\Id+s\xi)_{\#} \mu)-\Var(\mu)}{s}\right] \geq \INTDom{\langle x_1 - \bar{\mu} , \xi(x_1) \rangle}{\R^d}{\mu(x_1)}.
\end{equation*}

Furthermore, one can check that it holds
\begin{equation*}
\begin{aligned}
& \lim\limits_{s \downarrow 0} \left[ \frac{\Var((\Id+s\xi)_{\#} \mu)-\Var(\mu)}{s}\right] \\ 
\leq & \limsup\limits_{s \downarrow 0} \left[ \frac{\left(\Var((\Id+s\xi)_{\#} \mu)-\Var(\mu) \right)^+}{W_2(\mu,(\Id + s\xi)_{\#}\mu)} \right] \limsup\limits_{s \downarrow 0} \left[ \frac{W_2(\mu,(\Id + s\xi)_{\#}\mu)}{s}\right] \\
\leq & |\partial \Var|(\mu) \NormL{\xi}{2}{\mu},
\end{aligned}
\end{equation*}
so that, for any $\xi \in C^{\infty}_c(\R^d)$ with $\NormL{\xi}{2}{\mu} \leq 1$, one has
\begin{equation*}
 \INTDom{\langle x_1 - \bar{\mu} , \xi(x_1) \rangle}{\R^d}{\mu(x_1)} \leq |\partial \Var|(\mu).
\end{equation*}

By applying a density argument for test functions in the space $L^2(\R^d,\R^d;\mu)$ and using the dual characterization of the $L^2$-norm of a functional, it finally holds that $\NormL{\Id - \bar{\mu}}{2}{\mu} \leq |\partial \Var|(\mu)$, which amounts to state by Theorem \ref{thm:subdifferential_localslope} that the strong subdifferential $x \in \supp(\mu) \mapsto x - \bar{\mu}$ is the minimal selection in the classical subdifferential $\partial \Var(\mu)$ of the variance functional at $\mu$.
\end{proof}

\begin{rmk}[Possible extensions]
The analysis carried out in the previous Proposition for the variance functional can be applied in a similar fashion to integral functionals of the form
\begin{equation*}
\left\{
\begin{aligned}
& \Wcal^k : \mu \in \Pcal(K) \mapsto \INTDom{W(x_1,...,x_k)}{\R^d}{(\mu \times ...\times \mu)(x_1,...,x_k)} \\
& \Vcal_m : \mu \in \Pcal(K) \mapsto \INTDom{V \left( x,\mathsmaller{\INTDom{m(y)}{\R^d}{\mu(y)}} \right)}{\R^d}{\mu(x)}
\end{aligned}
\right.
\end{equation*}
for any $k \geq 1$, $W \in C^1(\R^{d \times k},\R)$, $V \in C^2(\R^d \times \R^n,\R)$ and $m \in C^2(\R^d,\R^n)$ for some $n \geq 1$.
\end{rmk}

\begin{prop}[Subdifferential of a smooth convolution interaction] Let $K \subset \R^d$ be a compact set, $H(\cdot,\cdot) \in C^1(\R^{2d},\R^d)$ be a function with sublinear growth and consider the non-local velocity field $v[\cdot](\cdot) : (\mu,x) \in \Pcal_c(K) \times \R^d \mapsto \INTDom{H(x,y)}{\R^d}{\mu(y)}$. 

Then, $v[\cdot](\cdot,\cdot)$ satisfies \textnormal{\textbf{(F)}}, \textnormal{\textbf{(B)}} and \textnormal{\textbf{(D)}} and the first order variations $x \in \supp(\mu) \mapsto \D_x v[\mu](x)$ and $x \in \supp(\mu) \mapsto \INTDom{\BGamma_{x}^{\circ}(y)}{\R^d}{\mu(y)}$ can be computed explicitly as
\begin{equation*}
\left\{
\begin{aligned}
& \D_x v[\mu](x) = \INTDom{\D_x H(x,y)}{\R^d}{\mu(y)}, \\
& \INTDom{\BGamma_{x}^{\circ}(y)}{\R^d}{\mu(y)} = \INTDom{\D_y H(x,y)}{\R^d}{\mu(y)}.
\end{aligned}
\right.
\end{equation*}
where $\BGamma^{\circ}_x(\cdot)$ is defined as in Theorem \ref{thm:General_PMP}.
\end{prop}

\begin{proof}
The Lipschitz estimates and the regularity in the sense of Definition \ref{def:Regular} can be derived using Kantorovich duality and the results of Proposition \ref{prop:Subdifferential_Integral}. For the first order variations, apply a classical differentiation under the integral sign result for the first one and Proposition \ref{prop:Subdifferential_Integral} to the components $\mu \mapsto \INTDom{H^i(x,y)}{\R^d}{\mu(y)}$ for any fixed $x \in \supp(\mu)$ for the second one. 
\end{proof}

We summarize these results in the form of an overview of possible functions satisfying \textbf{(H)} in the following corollaries.

\begin{cor}[Example of terminal costs satisfying the hypotheses of Theorem \ref{thm:General_PMP}]
\hfill \\ If $\varphi : \Pcal(K) \mapsto \R$ is either a (suitable) power of a smooth integral functional or the variance functional, then it satisfies hypotheses \textnormal{\textbf{(C)}}, \textnormal{\textbf{(B)}} and \textnormal{\textbf{(D)}}.
\end{cor}

\begin{cor}[Example of running costs satisfying the hypotheses of Theorem \ref{thm:General_PMP}]
\hfill \\ Let $l : (x,v) \in \R^{2d} \mapsto l(x,v) \in \R$ be a $C^1$, function $K \subset \R^d$ be compact, and define the running cost $L : (\mu,\omega) \in \Pcal(K) \times U \mapsto \R$ by
\begin{equation*}
L(\mu,\omega) = \INTDom{l(x,\omega(x))}{\R^d}{\mu(x)}.
\end{equation*}

Then, $L(\cdot,\cdot)$ satisfies hypotheses \textnormal{\textbf{(L)}\textnormal{}, \textbf{(B)}} and \textnormal{\textbf{(D)}}. Moreover, the barycenter of the minimal selection in its extended subdifferential $\Bpartial^{\circ}_{\mu} L(\mu,\omega)$ is determined at any $\mu \in \Pcal(K)$ by
\begin{equation*}
\bar{\Bgamma}^{\circ}_L \equiv \nabla_{\mu}L(\mu,\omega) : x \in \supp(\mu) \mapsto \nabla_x l(x,\omega(x)) + \D_x \omega(x)^{\top} \nabla_v l(x,\omega(x)).
\end{equation*}
\end{cor}

\begin{proof}
The proof only involves elementary Lipschitz-type estimates and the use of Proposition \ref{prop:Subdifferential_Integral}.
\end{proof}

Notice again that it is possible to take any power $\alpha \geq 1$ of the previous cost and any power $\alpha > 0$ provided that the functional does not vanish along the optimal pair control-trajectory $(u^*(\cdot),\mu^*(\cdot)) \in \Lip([0,T],\Pcal_c(\R^d)) \times \U$.

The following result shows that functionals based on kernels are regular. They appear in several mean-field models for interaction, see e.g. \cite{bellomo,Bellomo2013,CPT,MFOC,golse,spohn,vlasov}.

\begin{cor}[Non-local vector field satisfying hypotheses \textbf{(H)}]
\hfill \\ If $v[\cdot](\cdot,\cdot) : \Pcal(K) \times [0,T] \times K \rightarrow \R^d$ is defined for any $(\mu,t,x)$ by
\begin{equation*}
v[\mu](t,x) = (H(t,\cdot) \star \mu(t))(x) + v_l(t,x),
\end{equation*}
for some sublinear interaction kernel $H \in L^{\infty}([0,T],C^1(\R^d,\R^d))$ and vector field $v_l(\cdot,\cdot)$ measurable in $t$ as well as sublinear and Lipschitz in $x$, then it satisfies hypotheses \textnormal{\textbf{(F)}}, \textnormal{\textbf{(B)}} and \textnormal{\textbf{(C)}}.
\end{cor}

When the compactified Hamiltonian of the system $\omega \in U \mapsto \H_c(t,\nu^*(t),\omega) $ is differentiable at $u^*(t)$ for $\Lcal^1$-almost every $t \in [0,T]$, the maximization condition can be rewritten as a usual first-order condition.

\begin{cor}[Differentiation of the Hamiltonian]
Suppose that $u^*(t) \in \textnormal{int}(U)$ for $\Lcal^1$-almost every $t \in [0,T]$ and that the compactified Hamiltonian $\omega \in U \mapsto \H_c(t,\nu^*(t),\omega)$ is Fr\'echet differentiable at $u^*(t)$ for $\Lcal^1$-almost every $t \in [0,T]$. Then $u^*(t)$ verifies:
\begin{equation}
\D_{\omega} \H_c(t,\nu^*(t),u^*(t)) \cdot v = 0
\end{equation}
for all $v \in U$ and for $\Lcal^1$-almost every $t \in [0,T]$.
\end{cor}

\vspace{0.5cm}

\textbf{acknowledgements} This work has been carried out in the framework of Archimède
Labex (ANR-11-LABX-0033) and of the A*MIDEX project (ANR-
11-IDEX-0001-02), funded by the "Investissements d'Avenir" French
Government program managed by the French National Research
Agency (ANR). \\
The authors also thank the reviewer for the several useful comments that he/she provided.

%%%%%%%%%%%%%%%%%%%%%%%%%%%%%%%%%%%%%%%%%%%%%%%%%%%%%%%%%%%%%%%%%%%%%%%%%%%%%%%%
%							      BIBLIOGRAPHY                                 %
%%%%%%%%%%%%%%%%%%%%%%%%%%%%%%%%%%%%%%%%%%%%%%%%%%%%%%%%%%%%%%%%%%%%%%%%%%%%%%%%

{\footnotesize
\bibliographystyle{plain}
\bibliography{../../ControlWassersteinBib}}

\end{document}